\documentclass[11pt, twoside, final]{amsart}
\usepackage{etoolbox,lastpage}
\usepackage{amsmath,amsthm,amscd,amsfonts,amssymb,enumerate}
\usepackage{graphicx}
\usepackage{color}
\usepackage[colorlinks]{hyperref}
\usepackage{lineno}

\usepackage{graphicx}
\usepackage{geometry}
\usepackage{amsmath,amssymb,amsfonts,latexsym}
\usepackage{eurosym}
\usepackage[table]{xcolor}
\usepackage{float}
\usepackage{graphics}
\usepackage{subcaption}
\usepackage{listings,xcolor}
\usepackage{algorithm}
\usepackage{algorithmicx}
\usepackage{algpseudocode}

\newcommand{\WRP}{\par\quad\qquad\qquad\(\hookrightarrow\)\enspace}
\allowdisplaybreaks

\newtheorem{theorem}{Theorem}[section]

\newtheorem{corollary}[theorem]{Corollary}
\theoremstyle{definition}
\newtheorem{definition}[theorem]{Definition}

\theoremstyle{remark}
\newtheorem{remark}[theorem]{Remark}
\numberwithin{equation}{section}

\begin{document}
\title[Derivation of Computational Formulas for certain class of finite sums]{Derivation of Computational Formulas for certain class of finite sums: Approach to Generating functions arising from $p$-adic integrals and special functions}

\author{Yilmaz Simsek$^*$}
\address[Yilmaz Simsek]{Department of Mathematics, Faculty of Science
University of Akdeniz TR-07058 Antalya-TURKEY}
\email{ysimsek@akdeniz.edu.tr }
\thanks{$^*$Corresponding author}

\maketitle

\begin{abstract}
The aim of this paper is to construct generating functions for
some families of special finite sums with the aid of the Newton-Mercator
series, hypergeometric series, and $p$-adic integral (the Volkenborn
integral). By using these generating functions, their functional equations,
and their partial derivative equations, many novel computational formulas
involving the special finite sums of (inverse) binomial coefficients, the
Bernoulli type polynomials and numbers, Euler polynomials and numbers, the
Stirling numbers, the (alternating) harmonic numbers, the Leibnitz
polynomials and others. Among these formulas, by considering a computational formula which computes the aforementioned certain class of finite sums with the aid of the Bernoulli numbers and the Stirling numbers of the first kind, we present a computation algorithm and we provide some of their special values. Morover, using the aforementioned special finite sums and combinatorial numbers, we give relations among multiple alternating zeta
functions, the Bernoulli polynomials of higher order and the Euler
polynomials of higher order. We also give decomposition of the multiple
Hurwitz zeta functions with the aid of finite sums. Relationships and
comparisons between the main results given in the article and previously
known results have been criticized. With the help of the results of this
paper, the solution of the problem that Charalambides \cite[Exercise 30, p.
273]{Charalambos} gave in his book was found and with the help of this
solution, we also find very new formulas. In addition, the solutions of some
of the problems we have raised in \cite{SimsekMTJPAM2020} are also given.	

\end{abstract}

\section{Introduction, definitions and motivation}

It is known that there are a large number of researchers have studied to
find computational formulas for finite sums and infinite sums. Because it is
often not easy to find computational formulas for any finite sum, involving
special functions, special numbers, special polynomials, sums of higher
powers of binomial coefficients. In order to find any computational formula for finite sums, still many new methods and techniques have been developed,
investigated in mathematics, and also in other applied sciences. We know
that finite sums and their computational formulas are special important mathematical tools most used by mathematicians, physicists, engineers and
other scientists. Applications of the generating functions for special
numbers and polynomials, and finite sums with their computational formulas
have also been given by many different methods (\textit{cf}. \cite{Apostol}-%
\cite{Zave}). With this motivation, by an approach to generating functions arising from $p$-adic integrals and special functions, our purpose and motivations in this paper are to develop a computational methodology by deriving computational formulas for certain class of finite sums. The provided computational methodology provides the researchers a variety of methods that they can use in different fields and many situations.

In this study, we will construct generating functions that include special
numbers and polynomials, and special finite sums. With the help of these
generating functions and their functional equation, some new computational
formulas will be given for these special finite sums. On the other hand, the
main motivation of this paper is to construct and investigate generating
functions, given by Theorem \ref{Theorem1} and Theorem \ref{Theorem2}, with
a great deal distinct applications for the following numbers $y(n,\lambda )$%
, represented with certain finite sum:%
\begin{equation}
y(n,\lambda )=\sum_{j=0}^{n}\frac{(-1)^{n}}{(j+1)\lambda ^{j+1}\left(
	\lambda -1\right) ^{n+1-j}}  \label{ynldef}
\end{equation}%
(\textit{cf.} \cite{SimsekMTJPAM2020}).

The numbers $y(n,\lambda )$ reveal from the following zeta type function $%
\mathcal{Z}_{1}(s;a,\lambda )$:%
\begin{equation}
\mathcal{Z}_{1}(s;a,\lambda )=\frac{\ln \lambda }{\left( \lambda -1\right)
	(\ln a)^{s}}Li_{s}\left( \frac{1}{\lambda -1}\right) +\frac{1}{(\ln a)^{s}}%
\sum\limits_{n=0}^{\infty }\frac{y(n,\lambda )}{(n+1)^{s}},  \label{1aGZ}
\end{equation}%
where $\lambda \in \mathbb{C}\setminus \left\{ 0,1\right\} $ ($|\frac{1}{%
	\lambda -1}|<1$; $\operatorname{Re}(s)>1$), $a\geq 1$, and $s\in \mathbb{C}$, (%
\textit{cf}. \cite{SimsekMTJPAM2020}, \cite{SimsekREVISTA}, \cite%
{SimsekBull2021}). The function $Li_{s}\left( \lambda \right) $ given on the
right-hand side of the equation (\ref{1aGZ}) is denoted the polylogarithm
function:%
\begin{equation*}
Li_{s}(z)=\sum\limits_{n=1}^{\infty }\frac{\lambda ^{n}}{n^{s}},
\end{equation*}%
(\textit{cf}. \cite{SrivatavaChoi}, \cite{ChoiMTJPAM}).

Our first goal is to construct the following generating functions for the
numbers $y\left( n,\lambda \right) $ by aid the Newton-Mercator series which
is the Taylor series for the logarithm function.

\begin{definition}
	Let $\lambda \in \mathbb{R}$ with $\lambda \neq 0,1$. Let $z\in \mathbb{C}$.
	The numbers $y\left( n,\lambda \right) $ are defined by the following
	generating function:%
	\begin{equation}
	G\left( z,\lambda \right) =\sum\limits_{n=0}^{\infty }\left( 1-\lambda
	\right) ^{n+2}y\left( n,\lambda \right) z^{n}.  \label{1aG}
	\end{equation}
\end{definition}

This paper provide many new formulas that include not only the numbers $%
y\left( n,\lambda \right) $ and special numbers and polynomials and special
finite sums, but also their generating functions. Among others, we list the
following some of theorems involving these novel formulas. The proofs of the
Theorems are given in detail in the following sections.

\begin{theorem}
	\label{Theorem1} Let $\lambda \in \mathbb{R}$ with $\lambda \neq 0,1$. Let $%
	z\in \mathbb{C}$ with $\left\vert \frac{\lambda -1}{\lambda }z\right\vert <1$%
	. Then we have%
	\begin{equation}
	G\left( z,\lambda \right) =\frac{\ln \left( 1-\frac{\lambda -1}{\lambda }%
		z\right) }{z(z-1)}.  \label{1aG1}
	\end{equation}
\end{theorem}

To show the power of the function $G\left( z,\lambda \right) $ and its
various applications, it will be examined in detail in the next sections.

Subsequent, we come up with an interesting theorem asserting that the
function $G\left( z,\lambda \right) $ can be expressed as a hypergeometric
function $_{2}F_{1}\left[ 
\begin{array}{c}
1,1 \\ 
2%
\end{array}%
;\frac{1-\lambda }{\lambda }z\right] $.

\begin{theorem}
	\label{Theorem2} Let $\lambda \in \mathbb{R}$ with $\lambda \neq 0,1$. Let $%
	z\in \mathbb{C}$ with $\left\vert \frac{\lambda -1}{\lambda }z\right\vert <1$%
	. Then we have%
	\begin{equation}
	G\left( z,\lambda \right) =\frac{\left( 1-\lambda \right) z}{\lambda \left(
		z-1\right) }\text{ }_{2}F_{1}\left[ 
	\begin{array}{c}
	1,1 \\ 
	2%
	\end{array}%
	;\frac{1-\lambda }{\lambda }z\right] .  \label{1aGH}
	\end{equation}
\end{theorem}

Applying Theorem \ref{Theorem1} and Theorem \ref{Theorem2}, we derive many
novel computational formulas and relations involving finite sums, special
numbers, and special polynomials.

In addition, we present several eliciting and considerable remarks on these
formulas and relations.

By using generating functions and their functional equations, we give some
properties of the numbers $y(n,\lambda )$. We show that the numbers $%
y(n,\lambda )$ are closely associated with the Bernoulli numbers, the Euler
numbers, the harmonic numbers, the alternating Harmonic numbers, the
Apostol-Bernoulli numbers, the Stirling numbers, the Leibnitz numbers, and
special finite sums.

Another important purpose of this paper is to provide solutions to some of
the open problems that have been raised by the author \cite{SimsekMTJPAM2020}
and the Exercise 30 given by Charalambides \cite[Exercise 30, p. 273]%
{Charalambos} with the aid of the numbers $y(n,\lambda )$ and their
generating functions.

Some of our main contributions, derived from Theorem \ref{Theorem1}, Theorem %
\ref{Theorem2}, and hypergeometric functions, are listed by the following
Theorems, among other results.

\begin{theorem}
	\label{Theorem3}%
	Let $m\in \mathbb{N}_0$. Then, we have%
	\begin{equation}
	\mathcal{B}_{m}\left( \lambda \right) =\sum\limits_{n=0}^{m}\left(
	n+1\right) !\lambda ^{n+1}y\left( n,\lambda \right) S_{2}(m,n+1),
	\label{Abn-1}
	\end{equation}%
	where $\mathcal{B}_{m}\left( \lambda \right) $ and $S_{2}(m,n)$ denote the
	Apostol-Bernoulli numbers and the Stirling numbers of the second kind,
	respectively.
\end{theorem}

By applying  to the equation (\ref{1aG1}), we derive the following a new
relation involving the harmonic numbers $H_{n}$ and the numbers $y\left(
n,\lambda \right) $.

\begin{theorem}
	\label{Theorem H} Let $n\in \mathbb{N}$. Then we have%
	\begin{eqnarray}
	H_{2n+2}-H_{n+1} &=&-\frac{1}{2\left( n+1\right) }+\left( n+1\right)
	\sum\limits_{k=0}^{2n-1}\frac{\lambda ^{k+2}}{2n-k}y\left( n,\lambda \right)
	\label{aHY} \\
	&&+\left( n+1\right) \left( \lambda -1\right) \sum\limits_{k=0}^{2n}\frac{%
		\lambda ^{k+1}}{2n+1-k}y\left( n,\lambda \right) .  \notag
	\end{eqnarray}
\end{theorem}

\begin{theorem}
	\label{Theorem 4} Let $n\in \mathbb{N}$. Then we have%
	\begin{equation}
	\sum_{j=0}^{n}\frac{\zeta _{E}^{(n+1-j)}\left( -m,n+2\right) }{\left(
		j+1\right) 2^{n+1-j}}=\sum_{j=0}^{n}\frac{E_{m}^{(n+1-j)}\left( n+2\right) }{%
		\left( j+1\right) 2^{n+1-j}},  \label{ah5z}
	\end{equation}%
	where $\zeta _{E}^{(d)}\left( s,x\right) $ denotes the multiple alternating
	Hurwitz zeta function (multiple Hurwitz-Euler eta function), which is given
	as fallows \textup{(\ref{ah5zE})}:%
	\begin{equation}
	\zeta _{E}^{(d)}\left( s,x\right) =2^{d}\sum_{v=0}^{\infty }(-1)^{v}\binom{%
		v+d-1}{v}\frac{1}{\left( x+v\right) ^{s}}  \label{ah5zE}
	\end{equation}%
	(\textit{cf}. \cite{ChoiSrivastavaTJM}, \cite{SrivatavaChoi}).
\end{theorem}

Note that proof of Theorem \ref{Theorem 4} will be presented in Section \ref{Section6}.

\begin{theorem}
	\label{Theorem 5}%
	Let $m\in \mathbb{N}_0$. Then, we have%
	\begin{equation}
	y\left( m,\lambda \right) =\sum_{v=0}^{m}\sum_{n=0}^{v}(-1)^{v-m}\frac{%
		\left( \lambda -1\right) ^{v-m-1}B_{n}S_{1}(v,n)}{\lambda ^{v+1}v!}.
	\label{1aGbs}
	\end{equation}
\end{theorem}

We give two different proof of Theorem \ref{Theorem 5}. The first proof
related to generating functions and functional equation method. The second
proof is associated with the $p$-adic integral method.

\begin{theorem}
	\label{Theorem 6} Let $f(\lambda )$ be an entire function and $|\lambda |<1$%
	. Then we have%
	\begin{equation*}
	\sum_{v=0}^{\infty }f(v)\lambda ^{v}+\sum_{m=1}^{\infty }\frac{f^{(m-1)}(0)}{%
		m!}\sum\limits_{n=0}^{m}\left( n+1\right) !\lambda ^{n+1}y\left( n,\lambda
	\right) S_{2}(m,n+1)=0.
	\end{equation*}
\end{theorem}

\subsection{Preliminaries}

Throughout the paper, we use the following notation and definitions.

Let $\mathbb{N}$,$\mathbb{\ 
	\mathbb{Z}
}$, $\mathbb{R}$, and $\mathbb{C}$ denote the set of natural numbers, the
set of integer numbers, the set of real numbers and the set of complex
numbers, respectively. $\mathbb{N}_{0}:=\mathbb{N}\cup \left\{ 0\right\} $.
For $z\in \mathbb{C}$ with $z=x+iy$ ($x,y\in \mathbb{R}$); $\operatorname{Re}(z)=x$
and $\operatorname{Im}(z)=y$ and also $\ln z$ denotes the principal branch of the
many-valued function $\Im (\ln z)$ with the imaginary part of $\ln z$
constrained by%
\begin{equation*}
-\pi <\operatorname{Im}(\ln z)\leq \pi .
\end{equation*}%
\begin{equation*}
0^{n}=\left\{ 
\begin{array}{cc}
1, & n=0 \\ 
0, & n\in \mathbb{N}.%
\end{array}%
\right.
\end{equation*}

The well-known generalized hypergeometric function is given by%
\begin{equation}
_{p}F_{q}\left[ 
\begin{array}{c}
\alpha _{1},...,\alpha _{p} \\ 
\beta _{1},...,\beta _{q}%
\end{array}%
;z\right] =\sum\limits_{m=0}^{\infty }\left( \frac{\prod\limits_{j=1}^{p}%
	\left( \alpha _{j}\right) ^{\overline{m}}}{\prod\limits_{j=1}^{q}\left(
	\beta _{j}\right) ^{\overline{m}}}\right) \frac{z^{m}}{m!},  \label{hyper}
\end{equation}%
where the above series converges for all $z$ if $p<q+1$, and for $\left\vert
z\right\vert <1$ if $p=q+1$. Assuming that all parameters have general
values, real or complex, except for the $\beta _{j}$, $j=1,2,...,q$ none of
which is equal to zero or a negative integer and also 
\begin{equation*}
\left( \lambda \right) ^{\overline{v}}=\prod\limits_{j=0}^{v-1}(\lambda +j),
\end{equation*}%
and \ $\left( \lambda \right) ^{\overline{0}}=1$ for $\lambda \neq 1$, where 
$v\in \mathbb{N}$, $\lambda \in \mathbb{C}$. For the generalized hypergeometric function and their applications, it is
also recommended to refer to the following resource (\textit{cf}. \cite%
{Koepf,simsekJMAA,RT}, \cite{Rao}; and references therein).

\begin{equation*}
\left( \lambda \right) ^{\underline{v}}=\prod\limits_{j=0}^{v-1}(\lambda -j),
\end{equation*}%
and $\left( \lambda \right) ^{\underline{0}}=1$.

The Bernoulli polynomials of higher order, $B_{n}^{(l)}\left( y\right)$, are defined
by%
\begin{equation}
F_{B}(u,y;l)=\left( \frac{u}{e^{u}-1}\right) ^{l}e^{yt}=\sum_{n=0}^{\infty
}B_{n}^{(l)}\left( y\right) \frac{u^{n}}{n!},  \label{ApostolBern}
\end{equation}%
(\textit{cf}. \cite{balMS,CAC,comtet,KimDahee,Koepf,KucukogluAADM2019,Ozdemir,5Riardon,Roman,Rota,simsekAADM,Simsekfilomat,simsekJMAA,SimsekMTJPM,SimsekMTJPAM2020,1Sofo,SrivatavaChoi}; and references therein).

When $y=0$ in (\ref{ApostolBern}), we have the Bernoulli numbers of order $k$%
\begin{equation*}
B_{n}^{(k)}=B_{n}^{(k)}\left( 0\right) ,
\end{equation*}%
and when $k=1$, we have the Bernoulli numbers%
\begin{equation*}
B_{n}=B_{n}^{(1)},
\end{equation*}%
(\textit{cf}. \cite{CAC,comtet,KimDahee,Koepf,KucukogluAADM2019,Ozdemir,5Riardon,Roman,Rota,simsekAADM,Simsekfilomat,simsekJMAA,SimsekMTJPM,SimsekMTJPAM2020,1Sofo,SrivatavaChoi}; and references therein).

The Euler polynomials of higher order, $E_{m}^{(l)}(y)$, are defined by

\begin{equation}
F_{E}(u,y;l)=\left( \frac{2}{e^{u}+1}\right) ^{l}e^{yu}=\sum_{m=0}^{\infty
}E_{m}^{(l)}\left( y\right) \frac{u^{m}}{m!}  \label{Aeuler}
\end{equation}

The\ harmonic numbers $H_{n}$ are defined by%
\begin{equation}
F_{1}(u)=\frac{\ln (1-u)}{u-1}=\sum\limits_{n=1}^{\infty }H_{n}u^{n},
\label{aH}
\end{equation}%
where $H_{0}=0$ and $\left\vert u\right\vert <1$ (\textit{cf}. \cite{comtet,5Riardon,1Sofo,SrivatavaChoi,ChoiJIA}).

A relation between the numbers $y\left( n,\frac{1}{2}\right) $ and $H_{n}$
is given as follows:%
\begin{equation}
y\left( n,\frac{1}{2}\right) =2^{n+2}\left( H_{\left[ \frac{n}{2}\right]
}-H_{n}+\frac{(-1)^{n+1}}{n+1}\right) ,  \label{91c}
\end{equation}%
(\textit{cf}. \cite{SimsekREVISTA}).

The alternating Harmonic numbers $\mathcal{H}_{n}$ are defined by%
\begin{equation}
F_{2}(u)=\frac{\ln (1+u)}{u-1}=\sum\limits_{n=1}^{\infty }\mathcal{H}%
_{n}u^{n},  \label{AlH}
\end{equation}%
where $\left\vert u\right\vert <1$ (\textit{cf}. \cite{comtet}, \cite{Guo}, 
\cite{SimsekREVISTA}, \cite{1Sofo}).

In \cite[Eq. (20)]{SimsekREVISTA}, we showed that the following formula%
\begin{equation}
y\left( n,\frac{1}{2}\right) =2^{n+2}\sum_{j=0}^{n}\frac{(-1)^{j+1}}{j+1}
\label{91b}
\end{equation}%
is related to the following well-known alternating harmonic numbers%
\begin{equation*}
\mathcal{H}_{n}=\sum_{j=1}^{n}\frac{(-1)^{j}}{j}=H_{\left[ \frac{n}{2}\right]
}-H_{n},
\end{equation*}%
(\textit{cf}. \cite{comtet}, \cite{Guo}, \cite[Eq. (1.5)]{1Sofo}, \cite[Eq.
(20)]{SimsekREVISTA}).

The Stirling numbers of the first kind $S_{1}\left( v,d\right) $ are defined%
\begin{equation}
F_{s1}(u,d)=\frac{\left( \ln (1+u)\right) ^{d}}{d!}=\sum_{v=0}^{\infty
}S_{1}\left(v,d\right) \frac{u^{v}}{v!}  \label{Sitirling1}
\end{equation}%
and%
\begin{equation}
\left( u\right) ^{\underline{d}}=\sum_{j=0}^{d}S_{1}\left( d,j\right) u^{j}
\label{Sitirling1a}
\end{equation}%
(\textit{cf}. \cite{CAC,comtet,KimDahee,Koepf,KucukogluAADM2019,Ozdemir,5Riardon,Roman,Rota,simsekAADM,Simsekfilomat,simsekJMAA,simsekRSCM,SimsekMTJPM,SimsekMTJPAM2020,1Sofo,SrivatavaChoi}; and references therein).

Combining (\ref{Sitirling1}) with the Lagrange inversion formula, a
computation formula of the Stirling numbers of the first kind is given by%
\begin{equation}
S_{1}\left( n,k\right) =\sum_{c=0}^{n-k}\sum_{j=0}^{c}(-1)^{j}\binom{c}{j}%
\binom{n+c-1}{k-1}\binom{2n-k}{n-k-c}\frac{j^{n-k+c}}{c!},  \label{s1C}
\end{equation}%
where $k=0,1,2,\ldots ,n$ and $n\in \mathbb{N}_{0}$ (\textit{cf}. \cite[Eq.
(8.21), p. 291]{Charalambos}).

The Stirling numbers of the second kind $S_{2}\left( v,d\right) $ are
defined by%
\begin{equation}
F_{s2}(u,k)=\frac{\left( e^{u}-1\right) ^{d}}{d!}=\sum_{v=0}^{\infty
}S_{2}\left( v,d\right) \frac{u^{v}}{v!}  \label{S2}
\end{equation}%
(\textit{cf}. \cite{CAC,comtet,KimDahee,Koepf,KucukogluAADM2019,Ozdemir,5Riardon,Roman,Rota,simsekAADM,Simsekfilomat,simsekJMAA,simsekRSCM,SimsekMTJPM,SimsekMTJPAM2020,1Sofo,SrivatavaChoi}; and references therein).

Using (\ref{S2}), a computation formula of the Stirling numbers of the
second kind is given by%
\begin{equation}
S_{2}\left( n,k\right) =\frac{1}{k!}\sum_{c=0}^{k}(-1)^{k-c}\binom{k}{c}c^{n}
\label{S2C}
\end{equation}%
where $k=0,1,2,\ldots ,n$ and $n\in \mathbb{N}_{0}$ (\textit{cf}. \cite[Eq.
(8.19), p. 289]{Charalambos}).

The Apostol-Bernoulli numbers $\mathcal{B}_{v}\left( \theta \right) $ are
defined by%
\begin{equation}
F_{A}(u,\theta )=\frac{u}{\theta e^{u}-1}=\sum_{v=0}^{\infty }\mathcal{B}%
_{v}\left( \theta \right) \frac{u^{v}}{v!}  \label{Abn}
\end{equation}%
(\textit{cf}. \cite{Apostol}).

Combining (\ref{Abn}) with (\ref{S2}), a computation formula of the
Apostol-Bernoulli numbers $\mathcal{B}_{v}\left( \theta \right) $ is given by%
\begin{equation}
\mathcal{B}_{v}\left( \theta \right) =\frac{n\theta }{\left( \theta
	-1\right) ^{n}}\sum_{c=0}^{n-1}(-1)^{c}c!\theta ^{c-1}\left( \theta
-1\right) ^{n-1-c}S_{2}(n-1,c)  \label{AbnC}
\end{equation}
(\textit{cf}. \cite{Apostol}).

The Leibnitz polynomials $L_{m}(x)$ are defined by the following generating
function%
\begin{equation}
\mathcal{G}_{l}\left( x,u\right) =\frac{\ln \left( 1-u\right) +\ln \left(
	1-xu\right) }{\left( 1-u\right) \left( 1-xu\right) -1}=\sum\limits_{m=0}^{%
	\infty }L_{n}(x)u^{m},  \label{A1}
\end{equation}%
where $|u|<1$ (\textit{cf}. \cite[Exercise 16, p. 127]{Charalambos}).

Using (\ref{A1}), the polynomials $L_{m}(x)$, whose degree is $m$, are given
by%
\begin{equation*}
L_{m}(x):=\sum\limits_{l=0}^{m}\boldsymbol{l}\left( m,l\right) x^{l},
\end{equation*}%
where $\boldsymbol{l}\left( m,l\right) $\ denotes the Leibnitz numbers,
defined by%
\begin{equation}
\boldsymbol{l}\left( m,l\right) =\frac{1}{\left( m+1\right) \binom{m}{l}}
\label{ExpLeib}
\end{equation}%
or%
\begin{equation}
\boldsymbol{l}\left( m,l\right) =\sum\limits_{d=0}^{l}\left( -1\right) ^{l-d}%
\frac{1}{m-d+1}\binom{l}{d},  \label{SumLeib}
\end{equation}%
where $l=0,1,2,\ldots ,m$ and $m\in \mathbb{N}_{0}$ (\textit{cf}. \cite[%
Exercise 16, p. 127]{Charalambos}, \cite{SimsekASCM2021}).

The Bernstein basis functions are defined by%
\begin{equation}
B_{l}^{m}(x)=\binom{m}{l}x^{l}(1-x)^{m-l}  \label{aberns}
\end{equation}

Integrate the following equation (\ref{aberns}) with respect to $x$ from $0$
to $1$, we have%
\begin{equation*}
\frac{1}{\binom{m}{l}}\int\limits_{0}^{1}B_{l}^{m}(x)dx=\sum%
\limits_{d=0}^{m-l}(-1)^{m-l-d}\binom{m-l}{d}\int\limits_{0}^{1}x^{m-d}dx.
\end{equation*}%
Combining the above equation with (\ref{ExpLeib}), we have%
\begin{equation*}
\boldsymbol{l}\left( m,l\right) =\sum\limits_{d=0}^{m-l}(-1)^{m-l-d}\binom{%
	m-l}{d}\frac{1}{m-d+1}.
\end{equation*}

The numbers $Y_{v}(\lambda )$ are defined by%
\begin{equation}
g(u;\lambda )=\frac{2}{\theta ^{2}u+\theta -1}=\sum_{v=0}^{\infty
}Y_{v}(\theta )\frac{u^{v}}{v!},  \label{1aYY}
\end{equation}%
where%
\begin{equation*}
Y_{v}(\lambda )=-\frac{2\left( v!\right) \theta ^{2v}}{\left( 1-\theta
	\right) ^{v+1}}
\end{equation*}%
(\textit{cf}. \cite[Eq. (2.13)]{SimsekTJM2018}).

The Bernoulli numbers of the second kind (the Cauchy numbers) $b_{v}(0)$ are
defined by%
\begin{equation}
F_{b2}(u)=\frac{u}{\ln (1+u)}=\sum_{n=0}^{\infty }b_{v}(0)\frac{u^{v}}{v!},
\label{Be-1t}
\end{equation}%
(\textit{cf}. \cite{comtet}, \cite[p. 116]{Roman}, \cite{SimsekMTJPM}, \cite%
{SrivatavaChoi}; see also the references cited in each of these earlier
works).

The numbers $D_{n}$, which are so-called the Daehee numbers, are defined by%
\begin{equation}
F_{3}(u)=\frac{\ln (1+u)}{u}=\sum_{n=0}^{\infty }D_{n}\frac{u^{n}}{n!},\quad\text{
}(u\neq 0,\left\vert u\right\vert <1)  \label{Da}
\end{equation}%
(\textit{cf}. \cite{KimDahee}).

By combining the Newton-Mercator series with (\ref{Da}), one has the
following forula:%
\begin{equation}
D_{n}=(-1)^{n}\frac{n!}{n+1},  \label{D}
\end{equation}%
(\textit{cf}. \cite[p. 117]{CAC}, \cite{KimDahee}, \cite[p. 45, Exercise 19
(b)]{5Riardon}).

The derangement numbers $d_{m}$ are defined by the following generating function:%
\begin{equation}
F_{d}(u)=\frac{e^{-u}}{1-u}=\sum_{m=0}^{\infty }d_{m}\frac{u^{m}}{m!},\quad\text{ 
}(\left\vert u\right\vert <1)  \label{dn}
\end{equation}
where%
\begin{equation*}
d_{m}=\sum\limits_{j=0}^{m}(-1)^{j}\left( m-j\right) !\binom{m}{j}
\label{dm}
\end{equation*}%
(\textit{cf}. \cite{Carlitz}, \cite{Ma}).

The Fibonacci-type polynomials in two variables are defined by the following
generating function:%
\begin{equation}
H\left( t;x,y;k,m,l\right) =\sum_{n=0}^{\infty }\mathcal{G}_{n}\left(
x,y;k,m,l\right) t^{n}=\frac{1}{1-x^{k}t-y^{m}t^{m+l}},  \label{GH}
\end{equation}%
where $k,m,l\in \mathbb{N}_{0}$ (\textit{cf}. \cite{ozdemirFilomat}).

Using (\ref{GH}), we have the following explicit formula for the polynomials 
$\mathcal{G}_{n}\left( x,y;k,m,l\right) $:%
\begin{equation*}
\mathcal{G}_{n}\left( x,y;k,m,l\right) =\sum_{c=0}^{\left[ \frac{n}{m+l}%
	\right] }\binom{n-c\left( m+l-1\right) }{c}y^{mc}x^{nk-mck-lck},
\end{equation*}%
where $\left[ a\right] $ is the largest integer $\leq a$ (\textit{cf}. \cite%
{ozdemirFilomat}, \cite{Ozdemir}).

Let's briefly summarize the next sections of the article.

In Section \ref{Section2}, we give the solution of the open problem 1, which has been
proposed by the author \cite[p.57, Open problem 1]{SimsekMTJPAM2020} about
the generating functions of the numbers $y\left( n,\lambda \right) $. We give many properties of this function.

In Section \ref{Section3}, with the help of generating functions and their functional equations, we
give many identities involving the
numbers $y(n,\lambda )$, the Bernoulli numbers of the second kind, the
harmonic numbers, alternating Harmonic numbers, the Apostol-Bernoulli
numbers, the Stirling numbers, the Leibnitz numbers, the Bernoulli numbers,
and sums involving higher powers of inverses of binomial coefficients.

In Section \ref{Section4}, we give computation algorithm for the numbers $y\left(
n,\lambda \right) $. We also give some values of the numbers $y\left(
n,\lambda \right) $.

In Section \ref{Section5}, we give differential equations of the generating functions and
their applications. We give some applications of these equation.

In Section \ref{Section6}, with the aid of the numbers $y\left( n,\lambda \right) $, we
give decomposition of the multiple Hurwitz zeta functions involving the
Bernoulli polynomials of higher order.

In Section \ref{Section7}, we give infinite series representations of the numbers $%
y\left( n,\lambda \right) $ on entire functions.

In Section \ref{Section8}, we also construct the generating function for the numbers $y\left( n,\lambda \right) $ with the help of Volkenborn integral on $p$-adic integers. We give some applications of the $p$-adic integral.

In Section \ref{Section9}, we conclusion about the results of this paper.

\section{Generating functions for the numbers $y\left( n,\protect\lambda %
	\right) $}
\label{Section2}

In this section, we construct generating functions for the numbers $y\left(
n,\lambda \right) $. We give some properties of these functions. By using
these functions and their functional equations, we give many new
computational formulas and relations for the numbers $y\left( n,\lambda
\right) $ and special finite sums.

We also give the solution of the following open problem, which has been
proposed in \cite[p.57, Open problem 1]{SimsekMTJPAM2020}:

\textit{What is generating function for the numbers }$y\left( n,2\right) $%
\textit{\ and the numbers }$y\left( n,\lambda \right) ?$

Its answer is given by Theorem \ref{Theorem1}.

\begin{proof}[Proof of Theorem \protect\ref{Theorem1}] Substituting (\ref{ynldef}) into (\ref{1aG}), after
	some calculations, we obtain%
	\begin{equation*}
	G\left( z,\lambda \right) =\sum\limits_{n=0}^{\infty }\sum_{j=0}^{n}\frac{1}{%
		j+1}\left( \frac{\lambda -1}{\lambda }\right) ^{j+1}z^{n}.
	\end{equation*}%
	The following result is obtained by decomposing the above series for which
	the Cauchy product of two infinite series has been applied.%
	\begin{equation*}
	G\left( z,\lambda \right) =\frac{1}{z}\sum_{n=0}^{\infty }\frac{(-1)^{n+1}}{%
		n+1}\left( \frac{1-\lambda }{\lambda }z\right)
	^{n+1}\sum\limits_{n=0}^{\infty }z^{n}.
	\end{equation*}%
	After combining the above equation with the following well-known
	Newton-Mercator series%
	\begin{equation*}
	\ln (1+z)=\sum_{j=0}^{\infty }\frac{(-1)^{j}}{j+1}z^{j+1},\text{ }%
	(\left\vert z\right\vert <1)
	\end{equation*}%
	and geometric series%
	\begin{equation*}
	\sum\limits_{n=0}^{\infty }z^{n}=\frac{1}{1-z},\text{ }(\left\vert
	z\right\vert <1)
	\end{equation*}%
	yields the assertation of Theorem \ref{Theorem1}.
\end{proof}

\begin{proof}[Proof of Theorem \protect\ref{Theorem2}]
	Substituting $p=2$, $q=1$, $\alpha _{1}=\alpha _{2}=1$, $\beta _{1}=1$ and $%
	u=\frac{1-\lambda }{\lambda }z$ into (\ref{hyper}), we obtain%
	\begin{equation*}
	_{2}F_{1}\left[ 
	\begin{array}{c}
	1,1 \\ 
	2%
	\end{array}%
	;\frac{1-\lambda }{\lambda }z\right] =\sum\limits_{n=0}^{\infty }\frac{%
		\left( 1\right) ^{\overline{m}}\left( 1\right) ^{\overline{m}}}{\left(
		2\right) ^{\overline{m}}m!}\left( \frac{1-\lambda }{\lambda }z\right) ^{m}.
	\end{equation*}%
	Multiplying both sides of the above equation by $\frac{1-\lambda }{\lambda
		\left( z-1\right) }$, we get%
	\begin{equation*}
	_{2}F_{1}\left[ 
	\begin{array}{c}
	1,1 \\ 
	2%
	\end{array}%
	;\frac{1-\lambda }{\lambda }z\right] \frac{1-\lambda }{\lambda \left(
		z-1\right) }=\frac{1}{z-1}\sum\limits_{n=0}^{\infty }\frac{1}{m+1}\left( 
	\frac{1-\lambda }{\lambda }z\right) ^{m+1}.
	\end{equation*}%
	Combining the above equation with (\ref{1aG1}), we obtain%
	\begin{equation*}
	\frac{\left( 1-\lambda \right) z}{\lambda \left( z-1\right) }\text{ }%
	_{2}F_{1}\left[ 
	\begin{array}{c}
	1,1 \\ 
	2%
	\end{array}%
	;\frac{1-\lambda }{\lambda }z\right] =G\left( z,\lambda \right) .
	\end{equation*}%
	Thus, proof of Theorem \ref{Theorem2} is completed.
\end{proof}

Some special cases of the generating function, given in (\ref{1aG1}), are
given as follows:

Substituting $\lambda =-1$ into (\ref{1aG1}), we have the following
generating function for the numbers $y\left( n,-1\right) $:%
\begin{equation}
g_{1}(z)=G\left( z,-1\right) =\frac{\ln (1-2z)}{z^{2}-z}=\sum\limits_{n=0}^{%
	\infty }2^{n+2}y\left( n,-1\right) z^{n},  \label{1aG2}
\end{equation}%
where $z\neq 0$, $z\neq 1$ and $\left\vert 2z\right\vert <1$.

The function $g_{1}(z)$ is associated with generating function for the
finite sums of powers of inverse binomial coefficients. These relationships
will be investigated in detail in the following Sections.

Substituting $\lambda =2$ into (\ref{1aG}), we have the following generating
function for the numbers $y\left( n,2\right) $:%
\begin{equation}
g_{2}(z)=G\left( z,2\right) =\frac{\ln (1-\frac{z}{2})}{z^{2}-z}%
=\sum\limits_{n=0}^{\infty }(-1)^{n}y\left( n,2\right) z^{n},  \label{1aG3}
\end{equation}%
where $z\neq 0$, $z\neq 1$ and $\left\vert z\right\vert <1$.

By using the function $g_{2}(z)$, relationships among the numbers $y\left(
n,2\right) $, the Bernoulli numbers, the Stirling numbers, and some special
finite sums will be investigated in detail in the following sections.

Substituting $\lambda =\frac{1}{2}$ into (\ref{1aG1}), we have the following
generating function for the numbers $y\left( n,\frac{1}{2}\right) $:%
\begin{equation}
g_{3}(z)=G\left( z,\frac{1}{2}\right) =\frac{\ln (1+z)}{z^{2}-z}%
=\sum\limits_{n=0}^{\infty }\frac{1}{2^{n+2}}y\left( n,\frac{1}{2}\right)
z^{n},  \label{1aYn}
\end{equation}%
where $z\neq 0$, $z\neq 1$ and $\left\vert z\right\vert <1$.

By using the function $g_{3}(z)$, relationships among the numbers $y\left( n,%
\frac{1}{2}\right) $, the Bernoulli numbers, the Stirling numbers, the
Harmonic numbers, and some special finite sums will be investigated in
detail in the following sections.

\section{Identities derived from Generating function}
\label{Section3}

In this section, using generating functions and their functional equations, we
give very interesting and novel formulas and identities involving the
numbers $y(n,\lambda )$, the Bernoulli numbers of the second kind, the
harmonic numbers, alternating Harmonic numbers, the Apostol-Bernoulli
numbers, the Stirling numbers, the Leibnitz numbers, the Bernoulli numbers,
and sums involving higher powers of inverses of binomial coefficients.

We begin this section by giving proofs of some theorems given in the Section 1.

\begin{proof}[Proof. of Theorem \protect\ref{Theorem H}]
	Multiply both sides of the equation (\ref{1aG}) by the function $\ln \left(
	1+\frac{\lambda -1}{\lambda }z\right) $ and, with the help of the
	Newton-Mercator series, we obtain%
	\begin{equation}
	\frac{\ln \left( 1-\frac{\lambda -1}{\lambda }z\right) \ln \left( 1+\frac{%
			\lambda -1}{\lambda }z\right) }{z(z-1)}=\sum\limits_{n=0}^{\infty
	}\sum\limits_{k=0}^{n}(-1)^{n}\frac{\left( 1-\lambda \right) ^{n+3}y\left(
		k,\lambda \right) }{\left( n+1-k\right) \lambda ^{n-k+1}}z^{n+1}.
	\label{Af-1}
	\end{equation}%
	By using Abel's summation formula, and using multiplication two of the
	Newton-Mercator series, Furdui \cite{Furdui} gave the following formula:%
	\begin{equation}
	\ln \left( 1-y\right) \ln \left( 1+y\right) =\sum\limits_{v=1}^{\infty
	}\left( H_{v}-H_{v}-\frac{1}{2v}\right) \frac{y^{2v}}{v},\text{ }y^{2}<1.
	\label{Af2}
	\end{equation}%
	Combining (\ref{Af-1}) with (\ref{Af2}), we get%
	\begin{eqnarray*}
		&&\sum\limits_{m=1}^{\infty }\left( H_{m}-H_{m}-\frac{1}{2m}\right) \left( 
		\frac{\lambda -1}{\lambda }\right) ^{2m}\frac{z^{2m}}{m} \\
		&=&\sum\limits_{n=3}^{\infty }\sum\limits_{k=0}^{n-3}(-1)^{n+1}\frac{\left(
			1-\lambda \right) ^{n}y\left( k,\lambda \right) }{\left( n-k-2\right)
			\lambda ^{n-k-2}}z^{n} \\
		&&+\sum\limits_{n=2}^{\infty }\sum\limits_{k=0}^{n-2}(-1)^{n+1}\frac{\left(
			1-\lambda \right) ^{n+1}y\left( k,\lambda \right) }{\left( n-k-1\right)
			\lambda ^{n-k-1}}z^{n}.
	\end{eqnarray*}%
	Therefore%
	\begin{eqnarray}
	&&\sum\limits_{n=1}^{\infty }\left( H_{n}-H_{n}-\frac{1}{2n}\right) \left( 
	\frac{\lambda -1}{\lambda }\right) ^{2n}\frac{z^{2n}}{n}  \label{1aGi} \\
	&=&\sum\limits_{n=1}^{\infty }\sum\limits_{k=0}^{2n-1}\frac{\left( 1-\lambda
		\right) ^{2n+2}y\left( k,\lambda \right) }{\left( 2n-k\right) \lambda ^{2n-k}%
	}z^{2n}+\sum\limits_{n=1}^{\infty }\sum\limits_{k=0}^{2n}\frac{\left(
		1-\lambda \right) ^{2n+3}y\left( k,\lambda \right) }{\left( 2n-k+1\right)
		\lambda ^{2n-k+}}z^{2n}  \notag \\
	&&\sum\limits_{n=1}^{\infty }\left( \sum\limits_{k=0}^{2n-2}\frac{\left(
		1-\lambda \right) ^{2n+1}y\left( k,\lambda \right) }{\left( 2n-k-1\right)
		\lambda ^{2n-k-1}}-\sum\limits_{k=0}^{2n-1}\frac{\left( 1-\lambda \right)
		^{2n+2}y\left( k,\lambda \right) }{\left( 2n-k\right) \lambda ^{2n-k}}%
	\right) z^{2n+1}.  \notag
	\end{eqnarray}%
	After making some necessary algebraic calculations in the previous equation,
	the coefficients of $z^{2n}$ are equalized and the above equation we arrive
	at the desired result.
\end{proof}

By using (\ref{1aGi}), we also get the following result:

\begin{corollary}
	Let $n\in\mathbb{N}$. Then we have
	\begin{equation*}
	\sum\limits_{k=0}^{2n-2}\frac{\left( 1-\lambda \right) ^{2n+1}y\left(
		k,\lambda \right) }{\left( 2n-k-1\right) \lambda ^{2n-k-1}}%
	-\sum\limits_{k=0}^{2n-1}\frac{\left( 1-\lambda \right) ^{2n+2}y\left(
		k,\lambda \right) }{\left( 2n-k\right) \lambda ^{2n-k}}=0.
	\end{equation*}
\end{corollary}

We give the following functional equation%
\begin{equation}
-\frac{\lambda -1}{\lambda }\left( x+1-\frac{\lambda -1}{\lambda }z\right) 
\mathcal{G}_{l}\left( x,\frac{\lambda -1}{\lambda }z\right) =(z-1)G\left(
z,\lambda \right) +x(xz-1)G\left( xz,\lambda \right) .  \label{A1l}
\end{equation}%
Combining (\ref{A1l}) with (\ref{1aG}) and (\ref{A1}) yields%
\begin{eqnarray*}
	&&-\left( x+1\right) \sum\limits_{n=0}^{\infty }L_{n}(x)\left( \frac{\lambda
		-1}{\lambda }\right) ^{n+1}z^{n}+x\sum\limits_{n=1}^{\infty
	}L_{n-1}(x)\left( \frac{\lambda -1}{\lambda }\right) ^{n+1}z^{n} \\
	&=&-\sum\limits_{n=0}^{\infty }(-1)^{n}\left( \lambda -1\right) ^{n+2}\left(
	x^{n+1}+1\right) y\left( n,\lambda \right) z^{n} \\
	&&-\sum\limits_{n=0}^{\infty }(-1)^{n}\left( 1-\lambda \right) ^{n+1}\left(
	x^{n+1}+1\right) y\left( n-1,\lambda \right) z^{n}.
\end{eqnarray*}%
Comparing the coefficients of $z^{n}$ on both sides of the above equation, we
have following theorem:
\begin{theorem}
	Let $n\in \mathbb{N}$. Then we have%
	\begin{equation}
	\left( x+1\right) L_{n}(x)-xL_{n-1}(x)=(-1)^{n}\lambda ^{n+1}\left(
	x^{n+1}+1\right) \left( \left( \lambda -1\right) y\left( n,\lambda \right)
	+y\left( n-1,\lambda \right) \right) .  \label{A1l1}
	\end{equation}
\end{theorem}

We set the following functional equation:%
\begin{equation*}
g_{3}(z)F_{b2}(z)=\frac{1}{z-1},
\end{equation*}%
where $\left\vert z\right\vert <1$. Combining the above functional equation
with (\ref{1aYn}) and (\ref{Be-1t}), we get%
\begin{equation*}
\sum\limits_{n=0}^{\infty }\frac{1}{2^{n+2}}y\left( n,\frac{1}{2}\right)
z^{n}\sum_{n=0}^{\infty }b_{n}(0)\frac{z^{n}}{n!}=-\sum_{n=0}^{\infty }z^{n}.
\end{equation*}%
Therefore%
\begin{equation*}
\sum\limits_{n=0}^{\infty }\sum_{j=0}^{n}\frac{1}{2^{j+2}}y\left( j,\frac{1}{%
	2}\right) \frac{b_{n-j}(0)}{(n-j)!}z^{n}=-\sum_{n=0}^{\infty }z^{n}.
\end{equation*}%
Comparing the coefficients of $z^{n}$ on both sides of the above equation, we
have following theorem:

\begin{theorem}
	Let $n\in \mathbb{N}_0$. Then, we have%
	\begin{equation*}
	\sum_{j=0}^{n}\frac{y\left( j,\frac{1}{2}\right) b_{n-j}(0)}{2^{j}(n-j)!}=-4.
	\end{equation*}
\end{theorem}

We give a decomposition of the generating function $g_{3}(z)$ as follows:
\begin{equation*}
g_{3}(z)=F_{2}(z)-F_{3}(z).
\end{equation*}%
Combining the above function with (\ref{Da}) and (\ref{AlH}), we obtain%
\begin{equation*}
\sum\limits_{n=0}^{\infty }\frac{1}{2^{n+2}}y\left( n,\frac{1}{2}\right)
z^{n}=\sum\limits_{n=1}^{\infty }\mathcal{H}_{n}z^{n}-\sum_{n=0}^{\infty
}D_{n}\frac{z^{n}}{n!}.
\end{equation*}%
Comparing the coefficients of $z^{n}$ on both sides of the above equation,
we have following theorem:
\begin{theorem}
	Let $n\in \mathbb{N}_0$. Then, we have%
	\begin{equation}
	y\left( n,\frac{1}{2}\right) =\frac{2^{n+2}}{n!}\left( n!\mathcal{H}%
	_{n}-D_{n}\right) .
	\end{equation}
\end{theorem}

\begin{corollary}
	Let $n\in \mathbb{N}_0$. Then, we have%
	\begin{equation}
	y\left( n,\frac{1}{2}\right) =\frac{2^{n+2}}{n!}\left( n!\mathcal{H}%
	_{n}-\sum_{j=0}^{n}B_{j}S_{1}(n,j)\right) .
	\end{equation}
\end{corollary}

\begin{corollary}
	Let $n\in \mathbb{N}_0$. Then, we have%
	\begin{equation}
	y\left( n,\frac{1}{2}\right) =2^{n+2}\left( \mathcal{H}_{n}+\frac{(-1)^{n+1}%
	}{n+1}\right) .  \label{1aYn1}
	\end{equation}
\end{corollary}

\begin{corollary}
	Let $n\in \mathbb{N}_0$. Then, we have%
	\begin{equation}
	y\left( n,\frac{1}{2}\right) =2^{n+1}\mathcal{H}_{n+1}.
	\end{equation}
\end{corollary}

We set%
\begin{equation*}
F_{2}(z)=zF(z).
\end{equation*}%
By using the above equation, we have%
\begin{equation*}
\sum\limits_{n=1}^{\infty }\mathcal{H}_{n}z^{n}=\sum\limits_{n=0}^{\infty }%
\frac{1}{2^{n+2}}y\left( n,\frac{1}{2}\right) z^{n+1}.
\end{equation*}%
After some elementary calculations, we arrive at the allowing result:

\begin{corollary}
	Let $n\in \mathbb{N}$. Then, we have%
	\begin{equation}
	\mathcal{H}_{n}=\frac{1}{2^{n+1}}y\left( n-1,\frac{1}{2}\right) .
	\label{aHhy}
	\end{equation}
\end{corollary}

Noting that with the aid of (\ref{1aYn1}), we give a series representation of
the function $F_{2}(z)-g_{3}(z)$ as follows:%
\begin{equation}
F_{2}(z)-g_{3}(z)=\sum\limits_{n=0}^{\infty }\frac{(-1)^{n}}{n+1}z^{n}.
\label{1aYn2}
\end{equation}

By combining (\ref{1aYn}) with (\ref{1aYn2}), we also arrive at (\ref{aHhy}).

It is time to give the first proof of Theorem \ref{Theorem 5}. The following proof is related to generating functions and functional equation method.

\begin{proof}[The first proof of Theorem \protect\ref{Theorem 5}]
	Substituting $l=1$, $u=\ln \left( 1+\frac{1-\lambda }{\lambda }z\right) $, $%
	y=0$, and $z\neq 0$\ into (\ref{ApostolBern}), after some elementary
	calculations, we obtain%
	\begin{equation}
	\frac{\ln \left( 1+\frac{1-\lambda }{\lambda }z\right) }{\frac{1-\lambda }{%
			\lambda }z}=\sum_{n=0}^{\infty }B_{n}\frac{(\ln \left( 1+\frac{1-\lambda }{%
			\lambda }z\right) )^{n}}{n!}.  \label{d8}
	\end{equation}%
	Combining the above equation with (\ref{1aG}) and (\ref{Sitirling1}), we get%
	\begin{equation*}
	\sum\limits_{m=0}^{\infty }\left( 1-\lambda \right) ^{m+2}y\left( m,\lambda
	\right) z^{m}=\frac{1}{z-1}\sum_{m=0}^{\infty
	}\sum_{n=0}^{m}B_{n}S_{1}(m,n)\left( \frac{1-\lambda }{\lambda }\right)
	^{m+1}\frac{z^{m}}{m!}
	\end{equation*}%
	since $S_{1}(m,n)=0$ if $n>m$. Assuming that $\left\vert z\right\vert <1$,
	we obtain%
	\begin{equation*}
	\sum\limits_{m=0}^{\infty }\left( 1-\lambda \right) ^{m+2}y\left( m,\lambda
	\right) z^{m}=-\sum_{v=0}^{\infty }z^{v}\sum_{m=0}^{\infty
	}\sum_{n=0}^{m}B_{n}S_{1}(m,n)\left( \frac{1-\lambda }{\lambda }\right)
	^{m+1}\frac{z^{m}}{m!}.
	\end{equation*}%
	After some elementary calculations, the above equation yields%
	\begin{equation*}
	\sum\limits_{m=0}^{\infty }\left( 1-\lambda \right) ^{m+2}y\left( m,\lambda
	\right) z^{m}=-\sum_{m=0}^{\infty }\sum_{v=0}^{m}\sum_{n=0}^{v}\frac{\left(1-\lambda\right) ^{v+1}}{\lambda ^{v+1}}\frac{B_{n}S_{1}(v,n)}{v!}z^{m}.
	\end{equation*}%
	Now equating the coefficients of $z^{m}$ on both sides of the above
	equation, we arrive at the desired result.
\end{proof}

Substituting $\lambda =\frac{1}{2}$ and $\lambda =2$ into (\ref{1aGbs}), we
arrive at the following corollaries, respectively:

\begin{corollary}
	Let $m\in \mathbb{N}_0$. Then, we have%
	\begin{equation}
	y\left( m,\frac{1}{2}\right) =-2^{m+2}\sum_{v=0}^{m}\sum_{n=0}^{v}\frac{%
		B_{n}S_{1}(v,n)}{v!}.  \label{1aGbs1}
	\end{equation}
\end{corollary}

\begin{corollary}
	Let $m\in \mathbb{N}_0$. Then, we have%
	\begin{equation*}
	y\left( m,2\right) =\sum_{v=0}^{m}\sum_{n=0}^{v}\frac{\left( -1\right)
		^{v-m}B_{n}S_{1}(v,n)}{2^{v+1}v!}.
	\end{equation*}
\end{corollary}

\begin{remark}
	By combining the following well-known identity:%
	\begin{equation}
	\sum_{n=0}^{m}B_{n}S_{1}(m,n)=\frac{(-1)^{m}m!}{m+1}  \label{a91}
	\end{equation}%
	(\textit{cf}. \cite[p. 117]{CAC}, \cite{KimDahee}, \cite[p. 45, Exercise 19
	(b)]{5Riardon}, \cite[Eq. (20)]{SimsekREVISTA}), with (\ref{1aGbs}), we
	arrive at (\ref{ynldef}). Noting that there are many other proofs of (\ref%
	{a91}). For example, Kim \cite{KimDahee} gave proof of (\ref{a91}) by using the $p$-adic
	invariant integral on the set of $p$-adic integers. Kim represented the equation
	(\ref{a91}) by the notation $D_{n}$, which are so-called the Daehee numbers.
	Riordan \cite[p. 45, Exercise 19 (b)]{5Riardon} represented the equation (\ref{a91}) by the notation $(b)_{n}$.
\end{remark}

With the aid of the equations (\ref{a91}), (\ref{AlH}), (\ref{Da}) and (\ref%
{D}), we get some interesting formulas involving the Bernoulli numbers, the
Stirling numbers of the first kind, the Daehee numbers, and the alternating
Harmonic numbers.

Combining (\ref{AlH}) with (\ref{Da}), we get%
\begin{equation*}
F_{3}(u)=\frac{u-1}{u}F_{2}(u)
\end{equation*}%
Using the above equation, we get%
\begin{equation*}
\sum_{n=0}^{\infty }D_{n}\frac{u^{n+1}}{n!}=\sum\limits_{n=1}^{\infty }%
\mathcal{H}_{n}u^{n+1}-\sum\limits_{n=1}^{\infty }\mathcal{H}_{n}u^{n}.
\end{equation*}%
Therefore%
\begin{equation*}
\sum_{n=1}^{\infty }D_{n-1}\frac{u^{n}}{\left( n-1\right) !}%
=\sum\limits_{n=2}^{\infty }\mathcal{H}_{n-1}u^{n}-\sum\limits_{n=1}^{\infty
}\mathcal{H}_{n}u^{n}.
\end{equation*}%
Comparing the coefficients of $u^{n}$ on both sides of the above equation,
we arrive at the following relation:%
\begin{equation*}
D_{n-1}=\left( n-1\right) !\left( \mathcal{H}_{n-1}-\mathcal{H}_{n}\right) .
\end{equation*}%
By the above equation and (\ref{D}), we see that%
\begin{equation*}
\mathcal{H}_{n-1}-\mathcal{H}_{n}=\frac{(-1)^{n-1}}{n}.
\end{equation*}

Combining (\ref{AlH}) with (\ref{Da}), we also have%
\begin{equation*}
F_{2}(u)=-\frac{u}{1-u}F_{3}(u)
\end{equation*}%
Using the above equation, we obtain%
\begin{equation*}
\sum\limits_{n=1}^{\infty }\mathcal{H}_{n}u^{n}=-u\sum\limits_{n=0}^{\infty
}u^{n}\sum_{n=0}^{\infty }D_{n}\frac{u^{n}}{n!}.
\end{equation*}%
Therefore%
\begin{equation*}
\sum\limits_{n=1}^{\infty }\mathcal{H}_{n}u^{n}=-u\sum\limits_{n=0}^{\infty
}\sum_{j=0}^{n}D_{j}\frac{u^{n}}{j!}.
\end{equation*}%
Comparing the coefficients of $u^{n}$ on both sides of the above equation,
we arrive at the following relation:%
\begin{equation}
\mathcal{H}_{n}=-\sum_{j=0}^{n-1}\frac{D_{j}}{j!}.
\label{HarmonicDaehee}
\end{equation}%
Substituting (\ref{D}) into (\ref{HarmonicDaehee}), and using (\ref{a91}), we
arrive at the following theorem:

\begin{theorem}
	Let $n\in \mathbb{N}$. Then, we have%
	\begin{equation*}
	\mathcal{H}_{n}=-\sum_{j=0}^{n-1}\sum_{v=0}^{j}\frac{B_{v}S_{1}(j,v)}{j!}.
	\end{equation*}
\end{theorem}

It is time to give the first proof of Theorem \ref{Theorem3}. The following proof is related to generating functions and functional equation method.

\begin{proof}[Proof of Theorem \protect\ref{Theorem3}]
	Putting $z=\frac{\lambda }{1-\lambda }\left( e^{w}-1\right) $ in (\ref{1aG})
	and combining with (\ref{1aG1}), (\ref{Abn}) and (\ref{S2}), we obtain%
	\begin{equation*}
	F_{A}(w,\lambda )=\sum\limits_{n=0}^{\infty }\left( n+1\right) !\lambda
	^{n+1}y\left( n,\lambda \right) F_{s2}(w,n+1).
	\end{equation*}%
	By using the above equation, we get%
	\begin{eqnarray*}
		\sum_{m=0}^{\infty }\mathcal{B}_{m}\left( \lambda \right) \frac{w^{m}}{m!}
		&=&\sum\limits_{n=0}^{\infty }\left( n+1\right) !\lambda ^{n+1}y\left(
		n,\lambda \right) \sum_{m=0}^{\infty }S_{2}(m,n+1)\frac{w^{m}}{m!} \\
		&=&\sum_{m=0}^{\infty }\sum\limits_{n=0}^{m}\left( n+1\right) !\lambda
		^{n+1}y\left( n,\lambda \right) S_{2}(m,n+1)\frac{w^{m}}{m!}.
	\end{eqnarray*}%
	such that we here use the fact that $S_{2}(m,n)=0$ if $n>m$. Equating the
	coefficients of $\frac{w^{m}}{m!}$ on both sides of the above equation, we
	get the desired result.
\end{proof}

Combining (\ref{Abn-1}) with (\ref{ynldef}), we arrive at the following
corollaries:

\begin{corollary}
	Let $m\in \mathbb{N}_0$. Then, we have%
	\begin{equation*}
	\mathcal{B}_{m}\left( \lambda \right) =\frac{1}{\lambda -1}%
	\sum\limits_{n=0}^{m}\sum_{j=0}^{n}\frac{(-1)^{n}\left( n+1\right) !}{j+1}%
	\left( \frac{\lambda }{\lambda -1}\right) ^{n-j}S_{2}(m,n+1).
	\end{equation*}
\end{corollary}

\begin{corollary}
	Let $m\in \mathbb{N}_0$. Then, we have%
	\begin{equation*}
	\mathcal{B}_{m}\left( \lambda \right) =\frac{1}{1-\lambda }%
	\sum\limits_{n=0}^{m}\sum_{j=0}^{n}\frac{n+2}{j+1}\left( \frac{\lambda }{%
		\lambda -1}\right) ^{n-j}D_{n+1}S_{2}(m,n+1).
	\end{equation*}
\end{corollary}

\begin{corollary}
	Let $m\in \mathbb{N}_0$. Then, we have%
	\begin{equation*}
	\mathcal{B}_{m}\left( \lambda \right) =\frac{1}{1-\lambda }%
	\sum\limits_{n=0}^{m}\sum_{j=0}^{n}\sum_{k=0}^{n+1}\frac{n+2}{j+1}\left( 
	\frac{\lambda }{\lambda -1}\right) ^{n-j}B_{k}S_{1}(n+1,k)S_{2}(m,n+1).
	\end{equation*}
\end{corollary}

Substituting $z=e^{t}-1$ into (\ref{1aYn}), we get%
\begin{equation*}
\frac{t}{\left( e^{t}-2\right) \left( e^{t}-1\right) }=\sum\limits_{n=0}^{%
	\infty }\frac{1}{2^{n+2}}y\left( n,\frac{1}{2}\right) \left( e^{t}-1\right)
^{n}.
\end{equation*}

Combining the above equation with (\ref{Abn}) and (\ref{S2}), we get%
\begin{equation*}
\frac{1}{2t}\sum_{m=0}^{\infty }\mathcal{B}_{m}\left( \frac{1}{2}\right) 
\frac{t^{m}}{m!}\sum_{m=0}^{\infty }B_{m}\frac{t^{m}}{m!}=\sum_{m=0}^{\infty
}\sum\limits_{n=0}^{m}\frac{n!}{2^{n+2}}y\left( n,\frac{1}{2}\right)
S_{2}\left( m,n\right) \frac{t^{m}}{m!}
\end{equation*}%
or%
\begin{equation*}
\frac{1}{2}\sum_{m=0}^{\infty }\mathcal{B}_{m}\left( \frac{1}{2}\right) 
\frac{t^{m}}{m!}=2\sum_{m=0}^{\infty }\sum\limits_{n=0}^{m}\frac{(n+1)!}{%
	2^{n+2}}y\left( n,\frac{1}{2}\right) S_{2}\left( m,n+1\right) \frac{t^{m}}{m!%
}.
\end{equation*}%
After making some necessary algebraic calculations in the previous
equations, the coefficients of $\frac{t^{m}}{m!}$ are equalized, we arrive
at the following theorems:

\begin{theorem}
	Let $m\in \mathbb{N}_0$. Then, we have%
	\begin{equation}
	\mathcal{B}_{m}\left( \frac{1}{2}\right) =\sum\limits_{n=0}^{m}\frac{(n+1)!}{%
		2^{n+1}}y\left( n,\frac{1}{2}\right) S_{2}\left( m,n+1\right) .
	\label{1AAeq}
	\end{equation}
\end{theorem}

\begin{theorem}
	Let $m\in \mathbb{N}_0$. Then, we have%
	\begin{equation*}
	\sum_{n=0}^{m}\binom{m}{n}\mathcal{B}_{n}\left( \frac{1}{2}\right)
	B_{m-n}=m\sum\limits_{n=0}^{m-1}\frac{(n+1)!}{2^{n+1}}y\left( n,\frac{1}{2}%
	\right) S_{2}\left( m-1,n+1\right) .
	\end{equation*}
\end{theorem}

Substituting $u=e^{t}-1$ into (\ref{AlH}), we get%
\begin{equation*}
\frac{1}{2}\sum\limits_{m=0}^{\infty }\mathcal{B}_{n}\left( \frac{1}{2}%
\right) \frac{t^{m}}{m!}=\sum\limits_{m=0}^{\infty }\sum_{n=0}^{m}\frac{n!}{%
	m!}\mathcal{H}_{n}S_{2}(m,n)t^{m}.
\end{equation*}%
Combining the above equation with (\ref{1AAeq}), we arrive at the following
theorem:

\begin{theorem}
	Let $m\in \mathbb{N}_0$. Then, we have%
	\begin{equation*}
	\sum_{n=0}^{m}n!\mathcal{H}_{n}S_{2}(m,n)=\sum\limits_{n=0}^{m}\frac{(n+1)!}{2^{n}}%
	y\left( n,\frac{1}{2}\right) S_{2}\left( m,n+1\right) .
	\end{equation*}
\end{theorem}

Combining (\ref{1aG1}), (\ref{Da}) and (\ref{dn}), we get the following
functional equation:%
\begin{equation*}
\frac{\lambda-1 }{\lambda}F_{3}\left( \frac{1-\lambda }{\lambda }u\right)
F_{d}(u)=e^{-u}G(u,\lambda ).
\end{equation*}

By using the above equation, we obtain%
\begin{equation*}
\frac{\lambda-1 }{\lambda}\sum_{m=0}^{\infty }d_{m}\frac{u^{m}}{m!}%
\sum_{n=0}^{\infty }\left( \frac{1-\lambda }{\lambda }\right) ^{n}D_{n}\frac{%
	u^{n}}{n!}=\sum\limits_{n=0}^{\infty }\sum\limits_{j=0}^{n}\frac{(-1)^{n-j}}{%
	(n-j)!}\left( 1-\lambda \right) ^{j+2}y\left( j,\lambda \right) u^{n}.
\end{equation*}%
Therefore%
\begin{equation*}
-\sum_{n=0}^{\infty }\sum\limits_{m=0}^{n}\binom{n}{m}d_{m}\left( \frac{%
	1-\lambda }{\lambda }\right) ^{n-m+1}D_{n-m}\frac{u^{n}}{n!}%
=\sum\limits_{n=0}^{\infty }\sum\limits_{j=0}^{n}\frac{(-1)^{n-j}}{(n-j)!}%
\left( 1-\lambda \right) ^{j+2}y\left( j,\lambda \right) u^{n}.
\end{equation*}%
Comparing the coefficients of $u^{n}$ on both sides of the above equation,
we have following theorem:
\begin{theorem}
	Let $n\in \mathbb{N}_0$. Then, we have%
	\begin{equation*}
	-\sum\limits_{m=0}^{n}\binom{n}{m}\left( \frac{1-\lambda }{\lambda }\right)
	^{n-m+1}D_{n-m}d_{m}=\sum\limits_{j=0}^{n}\frac{(-1)^{n-j}}{(n-j)!}\left(
	1-\lambda \right) ^{j+2}y\left( j,\lambda \right) .
	\end{equation*}
\end{theorem}

Using (\ref{D}) and (\ref{dm}) in the above equation, we have%
\begin{eqnarray*}
	&&\sum\limits_{j=0}^{n}\frac{(-1)^{n-j}}{(n-j)!}\left( 1-\lambda \right)
	^{j+2}y\left( j,\lambda \right)  \\
	&=&\sum\limits_{m=0}^{n}\sum\limits_{j=0}^{m}(-1)^{n-m+j+1}\binom{n}{m}%
	\binom{m}{j}\left( \frac{1-\lambda }{\lambda }\right) ^{n-m+1}\frac{\left(
		n-m\right) !\left( m-j\right) !}{n-m+1}.
\end{eqnarray*}%
After some elementary calculations, we also arrive at the following corollary:
\begin{corollary}
	Let $n\in \mathbb{N}_0$. Then, we have%
	\begin{eqnarray*}
		\sum\limits_{j=0}^{n}\frac{(-1)^{n-j}}{(n-j)!}\left( 1-\lambda \right)
		^{j+2}y\left( j,\lambda \right)  =\sum\limits_{m=0}^{n}\sum\limits_{j=0}^{m}(-1)^{n-m+j+1}\left( \frac{%
			1-\lambda }{\lambda }\right) ^{n-m+1}\frac{n!}{\left( n-m+1\right) j!}.
	\end{eqnarray*}
\end{corollary}

\section{Computation algorithm for the numbers $y\left( n,\protect\lambda \right)$}
\label{Section4}

In this section, with the aid of (\ref{1aGbs}), (\ref{s1C}), and the
definition of the Bernoulli numbers, we present a computation algorithm (Algorithm \ref{alg:yCalculation} with a procedure called \texttt{COMPUTE$\char`_$y$\char`_$NUMBER}) for the numbers $y\left(m,\lambda \right) $.

\begin{algorithm}[H]
	\caption{Let $m$ be nonnegative integer and $\lambda\in \mathbb{C}$. This algorithm includes a procedure called \texttt{COMPUTE$\char`_$y$\char`_$NUMBER} which returns the numbers $y\left( m,\lambda \right) $.}
	\label{alg:yCalculation}
	\begin{algorithmic}
		\Procedure{\texttt{\textbf{COMPUTE$\char`_$y$\char`_$NUMBER}}}{$m$: nonnegative integer, $\lambda$}
		\State {$\textbf{Local variables: } v,n,y$}
		\State {$v,n,y\leftarrow 0$} 		
		\ForAll {$v$ in $\{0,1,2,\dots,m\}$}
		\ForAll {$n$ in $\{0,1,2,\dots,v\}$}
		\State $y\leftarrow y+\bigg(\Big($\texttt{Power}$\left(-1, v-m\right)*$\texttt{Power}$\left(\lambda-1, v-m-1\right)*$\texttt{BERNOULLI\_NUM}$\left(n\right)$\WRP$*$\texttt{STIRLING\_NUM\_FIRST}$\left(v,n\right)\Big)/\Big($\texttt{Power}$\left(\lambda, v+1\right)*$\texttt{Factorial}$\left(v\right)\Big)\bigg)$
		\EndFor
		\EndFor
		\State \textbf{return} {$y$}
		\EndProcedure
	\end{algorithmic}
\end{algorithm}

\begin{remark}
	In Algorithm \ref{alg:yCalculation}, the procedure \texttt{BERNOULLI$\char`_$NUM$\left(n\right)$} corresponds to the procedure which gives the $n$-th Bernoulli number. In addition, the procedure \texttt{STIRLING$\char`_$FIRST$\char`_$NUM} is corresponding to the procedure which computes the Stirling numbers of the first kind using the formula given in (\ref{s1C}). For details about the procedure \texttt{STIRLING$\char`_$FIRST$\char`_$NUM}, the interested readers may refer to the paper \cite{KucukogluAADM2019}.
\end{remark}

By using the computation
algorithm (Algorithm \ref{alg:yCalculation}), we give some values of the numbers $y\left(m,\lambda \right) $ as follows:
\begin{eqnarray*}
	y\left(0,\lambda \right)&=&\frac{1}{\lambda\left(\lambda-1 \right)},\\
	y\left(1,\lambda \right)&=&\frac{-3\lambda+1}{2\lambda^2\left(\lambda-1 \right)^2},\\
	y\left(2,\lambda \right)&=&\frac{11\lambda^2 -7\lambda +2}{6\lambda^3\left(\lambda-1 \right)^3},\\
	y\left(3,\lambda \right)&=&\frac{-25\lambda^3 +23\lambda^2 -13\lambda +3}{12\lambda^4\left(\lambda-1 \right)^4},\\
	y\left(4,\lambda \right)&=&\frac{137\lambda^4 -163\lambda^3 + 137\lambda^2 -63\lambda +12}{60\lambda^5\left(\lambda-1 \right)^5},
\end{eqnarray*}
and so on.
Observe that $y\left( n,\lambda \right) $ is a rational function of the
variable $\lambda $. The sequence of the leading coefficients of the
polynomial in the numerator of the numbers $y\left( n,\lambda \right) $ is
given as follows:

\begin{equation*}
1,-3,11,-25,137,-147,1089,-2283,7129,-7381,83711,\ldots 
\end{equation*}%
and so on. Taking absulute value of the each term above sequence, we get the
following well-known sequence:%
\begin{equation*}
\left( a(n)\right) _{n=1}^{\infty }=\left\{
1,3,11,25,137,147,1089,2283,7129,7381,83711,\ldots \right\} 
\end{equation*}%
and so on. The sequence $a(n)$ is given by OEIS: A025529 with the following
explicit formula (\textit{cf}. \cite{OEIS}):
\begin{equation*}
a(n)=\operatorname{lcm}(1,2,3,\ldots ,n)H_{n},
\end{equation*}%
where $n\in \mathbb{N}$, $H_{n}$ denotes the harmonic numbers.

\section{Differential equations of the generating functions and their
	applications}
\label{Section5}

In this section, we give partial derivative equations of the generating
functions. By applying these equations, we give many identities and many
novel recurrence relations involving the numbers $y(n,\lambda )$, $%
\boldsymbol{l}\left( n,0\right) $, and the special finite sums of (inverse) binomial coefficients.

Differentiating equation (\ref{1aG1}) with respect to $z$, we obtain the
following partial derivative equation:%
\begin{equation}
\left( z^{2}-z\right) \frac{\partial }{\partial z}\left\{ G\left( z,\lambda
\right) \right\} +\left( 2z-1\right) G\left( z,\lambda \right) =\frac{%
	1-\lambda }{\lambda +(1-\lambda )z}.  \label{1aGd1}
\end{equation}%
Differentiating equations (\ref{1aG2}), (\ref{1aG3}), and (\ref{1aYn}) with
respect to $z$, we obtain the following partial derivative equations,
respectively:%
\begin{equation}
\left( z^{2}-z\right) \frac{d}{dz}\left\{ g_{1}(z)\right\} +\left(
2z-1\right) g_{1}(z)=\frac{2}{2z-1},  \label{1aGd2}
\end{equation}%
\begin{equation}
\left( z^{2}-z\right) \frac{d}{dz}\left\{ g_{2}(z)\right\} +\left(
2z-1\right) g_{2}(z)=\frac{1}{z-2},  \label{1aGd3}
\end{equation}%
and%
\begin{equation}
\left( z^{2}-z\right) \frac{d}{dz}\left\{ g_{3}(z)\right\} +\left(
2z-1\right) g_{3}(z)=\frac{1}{z+1}.  \label{1aGd4}
\end{equation}

\subsection{Recurrence relations derived from PDEs for the generating
	functions}

Here, using equations (\ref{1aGd1})-(\ref{1aGd4}), we give recurrence
relations and identities involving the numbers $y(n,\lambda )$, $\boldsymbol{%
	l}\left( n,0\right) $, and the special finite sums of (inverse) binomial
coefficients.

\begin{theorem}
	The numbers $y(n,\lambda )$ satisfy the following derivative equations:%
	\begin{equation}
	(\lambda -1)\frac{d}{d\lambda }\left\{ y(n,\lambda )\right\}
	+(n+2)y(n,\lambda )=(-1)^{n+1}\sum\limits_{k=0}^{n}\frac{(\lambda -1)^{k-n-1}%
	}{\lambda ^{k+2}}.  \label{ynldefQED}
	\end{equation}%
	and%
	\begin{equation*}
	\frac{d}{d\lambda }\left\{ y(n,\lambda )\right\} +\frac{n+2}{\lambda -1}%
	y(n,\lambda )=\frac{(-1)^{n}}{\lambda }\left( 1-\left( \frac{\lambda }{%
		\lambda -1}\right) ^{n+1}\right) .
	\end{equation*}
\end{theorem}

\begin{proof}
	Differentiating equation (\ref{1aG1}) with respect to $\lambda $, we obtain%
	\begin{eqnarray*}
		\sum\limits_{n=0}^{\infty }(1-\lambda )^{n+1}\left( (1-\lambda )^{n+1}%
		\frac{d}{d\lambda }\left\{ y(n,\lambda )\right\} -(n+2)y(n,\lambda )\right)
		z^{n} =\frac{\partial }{\partial \lambda }\left\{ \frac{\ln \left( 1-\frac{%
				\lambda -1}{\lambda }z\right) }{z(z-1)}\right\} .
	\end{eqnarray*}%
	Therefore%
	\begin{eqnarray*}
		\lambda ^{2}\sum\limits_{n=0}^{\infty }(1-\lambda )^{n+1}\left( (1-\lambda
		)^{n+1}\frac{d}{d\lambda }\left\{ y(n,\lambda )\right\} -(n+2)y(n,\lambda
		)\right) z^{n} =\frac{1}{\left( 1-\frac{\lambda -1}{\lambda }z\right) \left( 1-z\right) }.
	\end{eqnarray*}%
	Assuming that $|\frac{\lambda -1}{\lambda }z|<1$ and $|z|<1$, then we obtain%
	\begin{eqnarray*}
		\lambda ^{2}\sum\limits_{n=0}^{\infty }(1-\lambda )^{n+1}\left( (1-\lambda
		)^{n+1}\frac{d}{d\lambda }\left\{ y(n,\lambda )\right\} -(n+2)y(n,\lambda
		)\right) z^{n} =\sum\limits_{n=0}^{\infty }\sum\limits_{k=0}^{n}\frac{(\lambda -1)^{k}}{%
			\lambda ^{k}}z^{n}.
	\end{eqnarray*}%
	Comparing the coefficients of $z^{n}$ on both sides of the above equation,
	we arrive at the desired result.
\end{proof}

Combining (\ref{1aG}) with (\ref{1aGd1}), we get%
\begin{eqnarray*}
	&&\left( z^{2}-z\right) \sum\limits_{n=1}^{\infty }n\left( 1-\lambda \right)
	^{n+2}y\left( n,\lambda \right) z^{n-1}+\left( 2z-1\right)
	\sum\limits_{n=0}^{\infty }\left( 1-\lambda \right) ^{n+2}y\left( n,\lambda
	\right) z^{n} \\
	&=&\sum\limits_{n=0}^{\infty }(-1)^{n}\left( \frac{1-\lambda }{\lambda }%
	\right) ^{n+1}z^{n}.
\end{eqnarray*}%
After some calculations in the above equation, after that equating the
coefficients of $z^{n}$ on both sides of the final equation, we get%
\begin{equation}
y\left( n-1,\lambda \right) +(\lambda -1)y\left( n,\lambda \right) =\frac{%
	(-1)^{n}}{(n+1)\lambda ^{n+1}}.  \label{1aGa}
\end{equation}%
Combining the above equation with (\ref{a91}), we arrive at a recurrence relation for the numbers $y\left( n,\lambda \right) $ as in the following theorem:
\begin{theorem}
	Let $n\in \mathbb{N}$. Then we have%
	\begin{equation*}
	y\left( n-1,\lambda \right) +(\lambda -1)y\left( n,\lambda \right) =\frac{1}{%
		\lambda ^{n+1}n!}\sum_{j=0}^{n}B_{j}S_{1}(n,j).
	\end{equation*}
\end{theorem}

\begin{remark}
	Substituting $x=0$ into (\ref{A1l1}) and using the following well-known
	identity%
	\begin{equation*}
	L_{n}(0)=\boldsymbol{l}\left( n,0\right) =\frac{1}{n+1},
	\end{equation*}%
	we also arrive at the equation (\ref{1aGa}).
\end{remark}

\begin{remark}
	By using (\ref{1aG1}) and (\ref{1aG}), we get%
	\begin{equation*}
	\sum\limits_{n=0}^{\infty }\left( 1-\lambda \right) ^{n+2}y\left( n,\lambda
	\right) z^{n+1}-\sum\limits_{n=0}^{\infty }\left( 1-\lambda \right)
	^{n+2}y\left( n,\lambda \right) z^{n}=\sum\limits_{n=0}^{\infty }\left( 
	\frac{1-\lambda }{\lambda }\right) ^{n+1}\frac{z^{n}}{n+1}.
	\end{equation*}%
	Equating the coefficients of $z^{n}$ on both sides of the above equation, we
	also arrive at the equation (\ref{1aGa}).
\end{remark}

Substituting $\lambda =2$ into (\ref{1aGa}), we have%
\begin{equation}
y\left( n-1,2\right) +y\left( n,2\right) =\frac{(-1)^{n}}{\left( n+1\right)
	2^{n+1}}.  \label{A1l2}
\end{equation}
Combining the above equation with the following well-known identity%
\begin{equation*}
\frac{(-1)^{n}}{2^{n}}n!=\sum_{j=0}^{n}E_{j}S_{1}(n,j)
\end{equation*}%
(\textit{cf}. \cite{Dkim}), we arrive at the following result:

\begin{corollary}
	Let $n\in \mathbb{N}$. Then we have%
	\begin{equation*}
	y\left( n-1,2\right) +y\left( n,2\right) =\frac{1}{2}\sum_{j=0}^{n}\frac{%
		E_{j}S_{1}(n,j)}{\left( n+1\right) !}.
	\end{equation*}
\end{corollary}

Combining (\ref{1aGd4}) with (\ref{1aYn}), and assuming that $\left\vert
z\right\vert <1$, we get%
\begin{eqnarray*}
	\left( z^{2}-z\right) \sum\limits_{n=1}^{\infty }\frac{1}{2^{n+2}}y\left( n,%
	\frac{1}{2}\right) z^{n-1} &=&\sum\limits_{n=0}^{\infty
	}(-1)^{n}z^{n}+\left( 1-2z\right) \\
	&&\times \sum\limits_{n=0}^{\infty }\frac{1}{2^{n+2}}y\left( n,\frac{1}{2}%
	\right) z^{n}.
\end{eqnarray*}%
After some elementary calculations in the above equation, we obtain%
\begin{eqnarray*}
	&&\sum\limits_{n=1}^{\infty }\frac{1}{2^{n+2}}y\left( n,\frac{1}{2}\right)
	z^{n+1}-\sum\limits_{n=1}^{\infty }\frac{1}{2^{n+2}}y\left( n,\frac{1}{2}%
	\right) z^{n} \\
	&=&\sum\limits_{n=0}^{\infty }(-1)^{n}z^{n}+\sum\limits_{n=0}^{\infty }\frac{%
		1}{2^{n+2}}y\left( n,\frac{1}{2}\right) z^{n}-2\sum\limits_{n=0}^{\infty }%
	\frac{1}{2^{n+2}}y\left( n,\frac{1}{2}\right) z^{n+1}.
\end{eqnarray*}%
Now equating the coefficients of $z^{n}$ on both sides of the above
equation, we arrive at the following corollary:
\begin{corollary}
	Let $n\in \mathbb{N}$. Then we have%
	\begin{equation}
	2y\left( n-1,\frac{1}{2}\right) -y\left( n,\frac{1}{2}\right) =(-1)^{n}\frac{%
		2^{n+2}}{n+1}.  \label{1aGd4a}
	\end{equation}
\end{corollary}

\begin{remark}
	Substituting $\lambda =\frac{1}{2}$ into (\ref{1aGa}), we also arrive at the
	equation (\ref{1aGd4}) and (\ref{1aGd4a}). Substituting $\lambda =2$ into (%
	\ref{1aGa}), we also arrive the (\ref{A1l2}).
\end{remark}

Substituting $\lambda =-1$ into (\ref{1aGa}), we get the following corollary:

\begin{corollary}
	Let $n\in \mathbb{N}$. Then we have%
	\begin{equation}
	y\left( n-1,-1\right) -2y\left( n,-1\right) =-\frac{1}{n+1}.  \label{1aGd1a}
	\end{equation}
\end{corollary}

Combining (\ref{1aGd1a}) with the following well-known formula%
\begin{equation*}
y\left( n,-1\right) =\frac{1}{2(n+1)}\sum\limits_{j=0}^{n}\frac{1}{\binom{n%
	}{j}}
\end{equation*}%
(\textit{cf}. \cite[Eq. (6.5)]{SimsekMTJPAM2020}), we get the following
combinatorial sum:

\begin{corollary}
	Let $n\in \mathbb{N}$. Then, we have%
	\begin{equation}
	\sum\limits_{j=0}^{n-1}\frac{1}{\binom{n-1}{j}}=\frac{2n}{n+1}%
	\sum\limits_{j=0}^{n-1}\frac{1}{\binom{n}{j}}.  \label{1aGf1}
	\end{equation}
\end{corollary}

Multiplying both sides of the equation (\ref{1aGd1a}) by $2(n+1)$, we arrive at
the following result:

\begin{corollary}
	Let $n\in \mathbb{N}$. Then we have%
	\begin{equation*}
	2(n+1)y\left( n-1,-1\right) -4(n+1)y\left( n,-1\right) =-2.
	\end{equation*}
\end{corollary}

\section{Decomposition of the multiple Hurwitz zeta functions with the help
	of the numbers $y\left( n,\protect\lambda \right) $}
\label{Section6}

In \cite{SimsekREVISTA}, by the aid of the numbers $y\left( n,\lambda
\right) $, we gave decomposition of the multiple Hurwitz zeta functions in
terms of the Bernoulli polynomials of higher order. In this section, by
using the same method in \cite{SimsekREVISTA}, we give decomposition of the
multiple alternating Hurwitz zeta functions in terms of the Bernoulli
polynomials of higher order, the Euler numbers and polynomials of higher
order, and the Stirling numbers of the first kind. By combining these
decomposition relations, we derive some formulas involving these numbers and
polynomials.

It is time to give the proof of Theorem \ref{Theorem 4} as follows:

\begin{proof}[Proof of Theorem \protect\ref{Theorem 4}]
	Substituting $\lambda =-e^{-t}$ into (\ref{ynldef}), we get%
	\begin{equation}
	y\left( n,-\frac{1}{e^{t}}\right) =\sum_{j=0}^{n}\frac{e^{t(n+2)}}{%
		(j+1)\left( e^{t}+1\right) ^{n+1-j}}.  \label{aH3a}
	\end{equation}%
	Combining (\ref{aH3a}) with (\ref{Aeuler}), we get 
	\begin{equation}
	y\left( n,-\frac{1}{e^{t}}\right) =\sum_{m=0}^{\infty }\sum_{j=0}^{n}\frac{%
		E_{m}^{(n+1-j)}\left( n+2\right) }{(j+1)2^{n+1-j}}\frac{t^{m}}{m!}.
	\label{aH3}
	\end{equation}%
	Using (\ref{aH3a}), we also get%
	\begin{equation}
	y\left( n,-\frac{1}{e^{t}}\right) =\sum_{j=0}^{n}\frac{1}{(j+1)}%
	\sum_{v=0}^{\infty }(-1)^{v}\binom{v+n-j}{v}e^{t(v+n+2)},  \label{aH3S}
	\end{equation}%
	where $\left\vert e^{t}\right\vert <1$. Using Taylor series of $e^{tx}$ in (%
	\ref{aH3S}) yields%
	\begin{equation}
	y\left( n,-\frac{1}{e^{t}}\right) =\sum_{j=0}^{n}\frac{1}{j+1}%
	\sum_{v=0}^{\infty }\sum_{m=0}^{\infty }(-1)^{v}\binom{v+n-j}{v}(v+n+2)^{m}%
	\frac{t^{m}}{m!}.  \label{aH4}
	\end{equation}%
	After making the required calculations in (\ref{aH3}) and (\ref{aH4}),
	equating the coefficients of $\frac{t^{m}}{m!}$ on both sides of the above
	equation, we obtain%
	\begin{equation}
	\sum_{j=0}^{n}\frac{1}{j+1}\left( \sum_{v=0}^{\infty }(-1)^{v}\binom{v+n-j}{v%
	}(v+n+2)^{m}-\frac{E_{m}^{(n+1-j)}\left( n+2\right) }{2^{n+1-j}}\right) =0.
	\label{ah5}
	\end{equation}%
	Therefore, proof is completed.
\end{proof}

\begin{remark}
	The well-known multiple Hurwitz-Euler eta function (or the multiple
	alternating Hurwitz function), which is given by the equation (\ref{ah5zE}),
	can also represent as follows:%
	\begin{equation*}
	\zeta _{E}^{(d)}\left( s,x\right) =2^{d}\sum_{v_{1},v_{2},\ldots
		,v_{d}=0}^{\infty }\frac{(-1)^{v_{1}+v_{2}+\cdots +v_{d}}}{\left(
		x+v_{1}+v_{2}+\cdots +v_{d}\right) ^{s}},
	\end{equation*}%
	where $\operatorname{Re}(s)>0$, $d\in \mathbb{N}$ and $x>0$ (\textit{cf}. \cite%
	{Cangul}, \cite{ChoiSrivastavaTJM}, \cite{Ozden}, 
	\cite{SimsekJdea}, \cite{SrivatavaChoi}).
\end{remark}

Combining the equation (\ref{ah5}) with the\ equation (\ref{ah5zE}), we also
arrive at the following theorem.

\begin{theorem}
	Let $m,n\in \mathbb{N}_{0}$. Then we have%
	\begin{eqnarray}
	&&\sum_{j=0}^{n}\frac{1}{\left( j+1\right) 2^{n+1-j}}\zeta
	_{E}^{(n+1-j)}\left( -m,n+2\right)  \label{1AAe} \\
	&=&\sum_{j=0}^{n}\sum_{l=0}^{m}\sum_{\underset{l_{1}+l_{2}+\cdots
			+l_{n+1-j}=l}{l_{1},l_{2},\cdots ,l_{n+1-j}=0}}^{l}\binom{m}{l}\frac{\left(
		n+2\right) ^{m-l}l!E_{l_{1}}E_{l_{2}}\cdots E_{l_{n+1-j}}}{%
		l_{1}!l_{2}!\cdots l_{n+1-j}!\left( j+1\right) 2^{n+1-j}}.  \notag
	\end{eqnarray}
\end{theorem}

Putting $n=0$ in (\ref{ah5z}), we get
\begin{equation*}
E_{m+1}(2)=\zeta _{E}(-m,2),
\end{equation*}%
where%
\begin{equation*}
\zeta _{E}(s,x)=\zeta _{E}^{(1)}(s,x)=2\sum\limits_{n=0}^{\infty }\frac{%
	(-1)^{n}}{(n+x)^{s}}
\end{equation*}%
(\textit{cf}. \cite{Cangul}, \cite{ChoiSrivastavaTJM}, \cite%
{Ozden}, \cite{SimsekJdea}, \cite{SrivatavaChoi}).

Substituting $n=1$ into (\ref{ah5z}), we obtain%
\begin{equation*}
\zeta _{E}^{(2)}(-m,3)-\zeta _{E}(-m,3)=E_{m}^{(2)}\left( 3\right)
-E_{m}\left( 3\right) .
\end{equation*}%
Substituting $n=2$ into (\ref{ah5z}), we also obtain%
\begin{eqnarray*}
	3\zeta _{E}^{(3)}(-m,4)+3\zeta _{E}^{(2)}(-m,4)+8\zeta _{E}(-m,3)=3E_{m}^{(3)}\left( 4\right) +3E_{m}^{(2)}\left( 4\right) +8E_{m}\left(
	4\right) .
\end{eqnarray*}%
Substituting $n=3$ into (\ref{aHh}), we also obtain%
\begin{eqnarray*}
	&&15\zeta _{E}^{(5)}(-m,5)+15\zeta _{E}^{(4)}(-m,5)+10\zeta
	_{E}^{(3)}(-m,5)+10\zeta _{E}^{(2)}(-m,5)+24\zeta _{E}(-m,5) \\
	&=&15E_{m}^{(5)}\left( 5\right) +15E_{m}^{(4)}\left( 5\right)
	+10E_{m}^{(3)}\left( 5\right) +10E_{m}^{(2)}\left( 5\right) +24E_{m}\left(
	5\right) .
\end{eqnarray*}

Since%
\begin{equation*}
\binom{v+n-j}{v}=\binom{v+n-j}{n-j},
\end{equation*}%
by combining the following well-known identity%
\begin{equation*}
\binom{v+n-j}{n-j}=\frac{1}{(n-j)!}\sum_{c=0}^{n-j}\left\vert
S_{1}(n-j,c+1)\right\vert v^{c},
\end{equation*}%
(\textit{cf. \cite{Choi1a}, \cite{comtet}, \cite{SrivatavaChoi}})\textit{\ }%
with the equation (\ref{aH4}), we get%
\begin{equation*}
y\left( n,-\frac{1}{e^{t}}\right)
=\sum_{j=0}^{n}\sum_{c=0}^{n-j}\sum_{v=0}^{\infty }\sum_{m=0}^{\infty }\frac{%
	(-1)^{v}(v+n+2)^{m}v^{c}}{\left( j+1\right) \left( n-j\right) !}\left\vert
S_{1}(n-j,c+1)\right\vert \frac{t^{m}}{m!}.
\end{equation*}

In \cite{SimsekREVISTA}, we gave the following results involving the
multiple Hurwitz zeta functions and the Bernoulli polynomials of higher
order:%
\begin{equation}
y\left( n,\frac{1}{e^{t}}\right) =\sum_{j=0}^{n}\frac{(-1)^{n}}{(j+1)}\zeta
_{n+1-j}\left( -m,n+2\right) \frac{t^{m}}{m!},  \label{aH3a1}
\end{equation}%
where $\zeta _{d}\left( s,x\right) $ denotes the Hurwitz zeta functions, for 
$d\in \mathbb{N}$, which defined by%
\begin{equation*}
\zeta _{d}\left( s,x\right) =\sum_{v=0}^{\infty }\binom{v+d-1}{v}\frac{1}{%
	\left( x+v\right) ^{s}}=\sum_{v_{1}=0}^{\infty }\sum_{v_{2}=0}^{\infty
}\cdots \sum_{v_{d}=0}^{\infty }\frac{1}{\left( x+v_{1}+v_{2}+\cdots
	+v_{d}\right) ^{s}}
\end{equation*}%
where $\Re (s)>d$, when $d=1$, we have the Hurwitz zeta function 
\begin{equation*}
\zeta (s,x)=\zeta _{1}\left( s,x\right) =\sum_{v=0}^{\infty }\frac{1}{%
	(x+v)^{s}},
\end{equation*}%
(\textit{cf}.  \cite{Choi1a}, \cite{MSKim1a}, \cite{SimsekJNT}-\cite{Simsek11a}, \cite{SrivatavaChoi}).

It is clear that%
\begin{equation}
\zeta _{d}\left( -m,x\right) =\frac{(-1)^{d}m!B_{m+d}^{(d)}(x)}{(d+m)!}
\label{aH3a4a}
\end{equation}%
and%
\begin{equation}
\zeta \left( -m,x\right) =\zeta _{1}\left( -m,x\right) =-\frac{B_{m+1}(x)}{%
	m+1}  \label{aH3a4}
\end{equation}%
where $m\in \mathbb{N}_{0}$ (\textit{cf}.  \cite{Choi1a}, 
\cite{MSKim1a}, \cite{SimsekJNT}-\cite{Simsek11a}, \cite%
{SrivatavaChoi}).

For $m,n\in \mathbb{N}$, we \cite{SimsekREVISTA} also defined 
\begin{equation}
\sum_{j=0}^{n}\frac{1}{j+1}\left( (-1)^{n}\zeta _{n+1-j}(-m,n+2)+\frac{%
	(-1)^{j}B_{m+n+1-j}^{(n+1-j)}\left( n+2\right) }{\binom{m+n+1-j}{n+1-j}%
	(n+1-j)!}\right) =0.  \label{aHh}
\end{equation}%
Substituting $\lambda =-e^{-2t}$ into (\ref{ynldef}), we get%
\begin{equation}
y\left( n,-\frac{1}{e^{2t}}\right) =\sum_{j=0}^{n}\frac{(-1)^{j-1}}{(j+1)}%
\frac{e^{2t(n+2)}}{\left( e^{2t}-1\right) ^{n+1-j}}.  \label{aH3eb}
\end{equation}%
Combining (\ref{aH3eb}) with (\ref{ApostolBern}) and (\ref{Aeuler}), we get%
\begin{equation}
y\left( n,-\frac{1}{e^{2t}}\right) =\sum_{j=0}^{n}\frac{(-1)^{j-1}}{%
	(j+1)2^{n+1-j}}\sum_{m=0}^{\infty }\frac{(-1)^{j}B_{m+n+1-j}^{(n+1-j)}\left(
	n+2\right) }{\binom{m+n+1-j}{n+1-j}(n+1-j)!}\frac{t^{m}}{m!}.  \label{aHeb}
\end{equation}%
Combining (\ref{aH3eb}) with (\ref{ApostolBern}), we also get%
\begin{equation}
y\left( n,-\frac{1}{e^{2t}}\right) =\sum_{j=0}^{n}\sum_{m=0}^{\infty
}\sum_{c=0}^{m}\frac{(-1)^{j-1}\binom{m}{c}B_{c+n+1-j}^{(n+1-j)}\left(
	2n+4\right) E_{m-c}^{(n+1-j)}\left( n+2\right) }{\binom{c+n+1-j}{n+1-j}%
	(n+1-j)!(j+1)2^{n+1-j}}\frac{t^{m}}{m!}.  \label{aH3Be}
\end{equation}

Combining (\ref{aHeb}) with (\ref{aH3Be}), we arrive at the following
theorem:

\begin{theorem}
	Let $n\in \mathbb{N}_0$. Then, we have%
	\begin{eqnarray*}
		&&\sum_{j=0}^{n}\frac{(-1)^{j-1}}{(n+1-j)!(j+1)2^{n+1-j}}\frac{%
			B_{m+n+1-j}^{(n+1-j)}\left( n+2\right) }{\binom{m+n+1-j}{n+1-j}} \\
		&=&\sum_{j=0}^{n}\frac{(-1)^{j-1}}{(n+1-j)!(j+1)2^{n+1-j}}\sum_{c=0}^{m}%
		\frac{\binom{m}{c}B_{c+n+1-j}^{(n+1-j)}\left( 2n+4\right) E_{m-c}^{(n+1-j)}}{%
			\binom{c+n+1-j}{n+1-j}}.
	\end{eqnarray*}
\end{theorem}

\begin{remark}
	Many other decompositions are obtained by continuing as above. It is known
	that the decomposition of the multiple Hurwitz zeta function is given by
	different techniques and methods in the literature. In this paper, we do not
	focus on the other kinds decompositions.
\end{remark}

\section{Infinite series representations of the numbers $y\left( n,\protect%
	\lambda \right) $ on entire functions}
\label{Section7}

In this section, we give some formulas containing the numbers $y\left(
n,\lambda \right) $ with the help of power series of entire functions. In
order to give these formulas, we need the following infinite series
representation, which was given by Boyadzhiev \cite{Boyadv1}-\cite{Boyadv2}:%
\begin{equation*}
\sum_{m=0}^{\infty }\frac{f^{(m)}(0)}{m!}h(m)y^{m}=\sum_{m=0}^{\infty }\frac{%
	h^{(m)}(0)}{m!}\sum\limits_{j=0}^{m}S_{2}(m,j)y^{j}f^{(j)}(y),
\end{equation*}%
where $f$ and $h$ are appropriate functions.

For a large class of entire functions and $|\lambda |<1$, Boyadzhiev \cite%
{Boyadv2} gave the following novel formula:%
\begin{equation}
\sum_{m=0}^{\infty }h(m)\lambda ^{m}+\sum_{m=1}^{\infty }\frac{h^{(m-1)}(0)}{%
	m!}\mathcal{B}_{m}(\lambda )=0.  \label{Abn-1b}
\end{equation}

Combining (\ref{Abn-1}) with (\ref{Abn-1b}), we arrive at the following
corollary:

\begin{corollary}
	Let $h(\lambda )$ be an entire function and $|\lambda |<1$. Then we have%
	\begin{equation}
	\sum_{v=0}^{\infty }h(v)\lambda ^{v}=-\sum_{m=1}^{\infty }\frac{h^{(m-1)}(0)%
	}{m!}\sum\limits_{n=0}^{m}\left( n+1\right) !\lambda ^{n+1}y\left( n,\lambda
	\right) S_{2}(m,n+1).  \label{ps1}
	\end{equation}
\end{corollary}

Substituting $h(\lambda )=\cos \lambda $ into (\ref{ps1}), with the aid of
the Euler formula, for $\left\vert \lambda e^{\pm i}\right\vert <1$, we
obtain%
\begin{equation}
\sum_{v=0}^{\infty }\lambda ^{v}\cos (v)=\frac{1}{2}\sum_{v=0}^{\infty
}\left( \left( \lambda e^{i}\right) ^{v}+\left( \lambda e^{-i}\right)
^{v}\right) =\frac{1-\lambda \cos 1}{1-2\lambda \cos 1+\lambda ^{2}}.
\label{GH1}
\end{equation}%
Combining (\ref{GH1}) with (\ref{ps1}), we have%
\begin{eqnarray}
&&\sum_{m=1}^{\infty }\sum\limits_{n=0}^{2m-1}\frac{(-1)^{m+1}\left(
	n+1\right) !S_{2}(2m-1,n+1)}{2(2m-1)!}\lambda ^{n+1}y\left( n,\lambda \right)
\label{Gh2} \\
&=&\frac{1-\lambda \cos 1}{1-2\lambda \cos 1+\lambda ^{2}},  \notag
\end{eqnarray}%
where $\left\vert \lambda \right\vert <1$.

Combining (\ref{Gh2}) with (\ref{GH}), we arrive at the following theorem:

\begin{theorem}
	\begin{eqnarray*}
		&&\sum_{m=1}^{\infty }\sum\limits_{n=0}^{2m-1}\frac{(-1)^{m+1}\left(
			n+1\right) !S_{2}(2m-1,n+1)}{2(2m-1)!}\lambda ^{n+1}y\left( n,\lambda \right)
		\\
		&=&1+\sum_{n=1}^{\infty }\left( \mathcal{G}_{n}(2\cos 1,-1;1,1,1)-\mathcal{G}%
		_{n-1}(2\cos 1,-1;1,1,1)\cos 1\right) \lambda ^{n}.
	\end{eqnarray*}
\end{theorem}

Noting that a considerable attribute of this formula is depended on the
function $f$. For instance, if $f$ is a polynomial, then infinite series on
the right-hand side of the equation (\ref{ps1}) reduces to the finite sum.
On account of this, one can easily compute the value of the infinite series
on the left-hand side of this equation.

The Hurwitz-Lerch zeta function $\Phi (\lambda ,z,b)$ is defined by%
\begin{equation*}
\Phi (\lambda ,z,b)=\sum_{j=0}^{\infty }\frac{\lambda ^{j}}{\left(
	j+b\right) ^{z}},
\end{equation*}%
where $b\in \mathbb{C}\setminus \mathbb{Z}_{0}^{-}$;$\;\lambda ,z\in \mathbb{%
	C}$ when $\;\left\vert \lambda \right\vert <1$;$\;x>1\;$when$\;\left\vert
\lambda \right\vert =1$ (\textit{cf}. \cite{Apostol}, \cite{SrivatavaChoi}).

The function $\Phi (\lambda ,z,b)$ interpolates the Apostol-Bernoulli
polynomials at negative integers, that is%
\begin{equation}
\Phi (\lambda ,1-n,b)=-\frac{1}{n}\mathcal{B}_{n}\left( b;\lambda \right) ,
\label{abp}
\end{equation}%
where $n\in \mathbb{N}$, $\left\vert \lambda \right\vert <1$ (\textit{cf}. 
\cite{Apostol}, \cite{Boyadv1}, \cite{Boyadv2}, \cite{SrivatavaChoi}).

Since%
\begin{equation*}
\mathcal{B}_{0}\left( 0;\lambda \right) =0
\end{equation*}%
and%
\begin{equation*}
\lambda \mathcal{B}_{1}(1;\lambda )=1+\mathcal{B}_{1}(\lambda )
\end{equation*}%
and for $n\geq 2$,%
\begin{equation*}
\lambda \mathcal{B}_{n}(1;\lambda )=\mathcal{B}_{n}(\lambda ),
\end{equation*}%
(\textit{cf}. \cite{Apostol}), for $b=1$, the equation (\ref{abp}) reduces
to the following well-known formula:%
\begin{eqnarray}
\Phi (\lambda ,1-n,1) &=&-\frac{1}{n}\mathcal{B}_{n}\left( 1;\lambda \right)
\label{1AAep} \\
&=&-\frac{1}{n\lambda }\mathcal{B}_{n}\left( \lambda \right) ,  \notag
\end{eqnarray}%
where $n\in \mathbb{N}$ with $n\geq 2$.

Substituting $f(z)=z^{m}$ ($m\in \mathbb{N}$) into (\ref{ps1}), and using (%
\ref{Abn-1}), we arrive at the another proof of Theorem \ref{Theorem3}.

\section{Construction of the generating function $G\left( t,\protect\lambda %
	\right) $ with the help of Volkenborn integral on $p$-adic integers}
\label{Section8}

In this section, we present another construction of the generating function $%
G\left( t,\lambda \right) $ with the help of Volkenborn integral on $p$-adic
integers. We give $p$-adic integral representation of the function $G\left(
t,\lambda \right) $. We also give some applications of this integral
representation.

Let $\mathbb{Z}_{p}$ denote the set of $p$-adic integers. Also, let $C^{1}$ denote the set of continuous differentiable functions from $\mathbb{Z}_{p}$ to a field with a complete valuation. With the aid of the indefinite sum of a continuous function $f$ on $\mathbb{Z}_{p}$, the Volkenborn integral is
given by%
\begin{equation}
\int\limits_{\mathbb{Z}_{p}}f\left( x\right) d\mu _{1}\left( x\right) =%
\underset{N\rightarrow \infty }{\lim }\frac{1}{p^{N}}\sum_{x=0}^{p^{N}-1}f%
\left( x\right) ,  \label{M}
\end{equation}%
where $f\in C^{1}$ on $\mathbb{Z}_{p}\ $and $\mu _{1}\left( x\right) $
denotes the Haar distribution, which is given by%
\begin{equation*}
\mu _{1}\left( x+p^{N}\mathbb{Z}_{p}\right) =\mu _{1}\left( x\right) =\frac{1%
}{p^{N}}
\end{equation*}%
(\textit{cf}. \cite{T. Kim}, \cite{KIMjmaa2017}, \cite[Definition 55.1, p.167%
]{Schikof}, \cite{SimsekMTJPM}, \cite{Volkenborn}); see also the references
cited in each of these earlier works).

\begin{theorem}
	Let $\lambda\in \mathbb{Z}_p$. Then, we have%
	\begin{equation}
	G\left( t,\lambda \right) =\frac{1-\lambda }{\lambda\left(
		t-1\right) }\int\limits_{\mathbb{Z}_{p}}\left( 1+\frac{1-\lambda}{\lambda }%
	t\right) ^{x}d\mu _{1}\left( x\right) .  \label{pG}
	\end{equation}
\end{theorem}

\begin{proof} Substituting $f(x;t,\lambda )=\left( 1+\frac{1-\lambda}{\lambda }%
	t\right) ^{x}$ into the following integral equation, which is given in \cite[
	p.169]{Schikof}:
	\begin{equation*}
	\int\limits_{\mathbb{Z}_{p}}f\left( x+1;t,\lambda \right) d\mu _{1}\left(
	x\right) =\int\limits_{\mathbb{Z}_{p}}f\left( x;t,\lambda \right) d\mu
	_{1}\left( x\right) +\frac{d}{dx}f(x)\left\vert _{x=0}\right. ,
	\end{equation*}%
	and after some elementary computations, we obtain%
	\begin{equation*}
	G\left( t,\lambda \right) =\frac{1-\lambda }{\lambda\left(
		t-1\right) }\sum\limits_{n=0}^{\infty }(-1)^{n}\left(\frac{1-\lambda}{%
		\lambda }t\right) ^{n}\frac{1}{n+1}.
	\end{equation*}%
	Combining the above equation with the Newton-Mercator series, we arrive at
	the desired result.
\end{proof}

It is time to give the second proof of Theorem \ref{Theorem 5}. The following proof is associated with the $p$-adic integral method.

\begin{proof}[The second proof of Theorem \protect\ref{Theorem 5}]
	Using (\ref{pG}), we get%
	\begin{equation}
	G\left( t,\lambda \right) =\frac{1-\lambda }{\lambda\left(
		t-1\right) }\sum\limits_{n=0}^{\infty }(-1)^{n}\left( \frac{\lambda -1}{%
		\lambda }t\right) ^{n}\frac{1}{n!}\int\limits_{\mathbb{Z}_{p}}\left(
	x\right) ^{\underline{n}}\mu _{1}\left( x\right) .  \label{pG2}
	\end{equation}%
	Combining the above equation with (\ref{Sitirling1a}), we obtain%
	\begin{equation}
	G\left( t,\lambda \right) =\frac{1-\lambda }{\lambda\left(
		t-1\right) }\sum\limits_{n=0}^{\infty }(-1)^{n}\left( \frac{\lambda -1}{%
		\lambda }t\right) ^{n}\frac{1}{n!}\sum\limits_{j=0}^{n}S_{1}(n,j)\int%
	\limits_{\mathbb{Z}_{p}}x^{j}\mu _{1}\left( x\right) .  \label{pG1}
	\end{equation}%
	Combining (\ref{pG1}) with the following well-known formula%
	\begin{equation*}
	B_{j}=\int\limits_{\mathbb{Z}_{p}}x^{j}\mu _{1}\left( x\right) 
	\end{equation*}%
	(\textit{cf}. \cite[ p.171]{Schikof}), and (\ref{1aG}), we have%
	\begin{equation*}
	\sum\limits_{n=0}^{\infty }\left( 1-\lambda \right) ^{n+2}y\left( n,\lambda
	\right) t^{n}=\frac{1-\lambda }{\lambda\left( t-1\right) }%
	\sum\limits_{n=0}^{\infty }\sum\limits_{j=0}^{n}\frac{S_{1}(n,j)B_{j}}{n!}%
	\left( \frac{1-\lambda }{\lambda }\right) ^{n}t^{n}.
	\end{equation*}%
	Therefore%
	\begin{equation*}
	\sum\limits_{n=0}^{\infty }\left( 1-\lambda \right) ^{n+2}y\left( n,\lambda
	\right) t^{n}=\sum\limits_{n=0}^{\infty
	}\sum\limits_{d=0}^{n}\sum\limits_{j=0}^{d}(-1)^{d}\frac{S_{1}(d,j)B_{j}}{%
		d!}\left( \frac{\lambda -1}{\lambda }\right) ^{d+1}t^{n}.
	\end{equation*}%
	Comparing the coefficients of $t^{n}$ on both sides of the above equation,
	we arrive at the equation (\ref{1aGbs}).
\end{proof}

By using (\ref{pG2}), we obtain%
\begin{eqnarray}
&&\sum\limits_{n=0}^{\infty }\left( 1-\lambda \right) ^{n+2}y\left(
n,\lambda \right) t^{n+1}-\sum\limits_{n=0}^{\infty }\left( 1-\lambda
\right) ^{n+2}y\left( n,\lambda \right) t^{n}  \label{pG3} \\
&=&\sum\limits_{n=0}^{\infty }(-1)^{n+1}\left( \frac{\lambda -1}{\lambda }%
\right) ^{n+1}\frac{1}{n!}\int\limits_{\mathbb{Z}_{p}}\left( x\right) ^{%
	\underline{n}}\mu _{1}\left( x\right) t^{n}.  \notag
\end{eqnarray}%
Combining the above equation with the Volkenborn integral in terms of the
Mahler coefficients%
\begin{equation*}
\int\limits_{\mathbb{Z}_{p}}\binom{x}{n}\mu _{1}\left( x\right) =\frac{%
	(-1)^{n}}{n+1}
\end{equation*}%
(\textit{cf}. \cite[Proposition 55.3, p.168]{Schikof}), we get%
\begin{eqnarray*}
	\sum\limits_{n=0}^{\infty }\left( 1-\lambda \right) ^{n+2}y\left(
	n,\lambda \right) t^{n+1}-\sum\limits_{n=0}^{\infty }\left( 1-\lambda
	\right) ^{n+2}y\left( n,\lambda \right) t^{n} =\sum\limits_{n=0}^{\infty }\frac{(-1)^{n}}{n+1}\left( \frac{1-\lambda }{%
		\lambda }\right) ^{n+1}t^{n}.
\end{eqnarray*}
Comparing the coefficients of $t^{n}$ on both sides of the above equation,
we arrive at the equation (\ref{1aGa}).

Combining the following well-known identity which proved by Kim et al. \cite%
{DSkimDaehee}:
\begin{equation*}
D_{m}=\int\limits_{\mathbb{Z}_{p}}\left( x\right) ^{\underline{m}}\mu
_{1}\left( x\right)
\end{equation*}%
where $m\in \mathbb{N}_{0}$, with (\ref{pG3}), we
get the following result:

\begin{theorem}
	Let $n\in \mathbb{N}$. Then, we have
	\begin{equation*}
	y\left( n-1,\lambda \right) +(\lambda -1)y\left( n,\lambda \right) =\frac{%
		D_{n}}{\lambda ^{n+1}n!}.
	\end{equation*}
\end{theorem}

\section{Conclusion}
\label{Section9}

In this paper, we have given the solution of the Open problem 1, which has been proposed by the author \cite[p.57, Open problem 1]{SimsekMTJPAM2020} about
the generating functions of the numbers $y\left( n,\lambda \right) $. We have also given many properties of this function.  with the help of generating functions and their functional equations, we
have derived many formulas associated with the
numbers $y(n,\lambda )$, the Bernoulli numbers of the second kind, the harmonic numbers, alternating Harmonic numbers, the Apostol-Bernoulli
numbers, the Stirling numbers, the Leibnitz numbers, the Bernoulli numbers,
and sums involving higher powers of inverses of binomial coefficients. Furthermore, we have provided an algorithm to compute the numbers $y\left(
n,\lambda \right) $. By using this algorithm, we have also computed some values of the numbers $y\left(
n,\lambda \right) $. In addition, we  have presented differential equations of the generating functions with their applications. With the aid of the numbers $y\left( n,\lambda \right) $, we
have given decomposition of the multiple Hurwitz zeta functions involving the Bernoulli polynomials of higher order. We have also given infinite series representations of the numbers $y\left( n,\lambda \right) $ on entire functions. By the aid of the Volkenborn integral on $p$-adic
integers, we have constructed the generating function for the numbers $y\left( n,\lambda \right) $ and given their some applications.

As a result, the results produced in this article are in a wide range and have the potential to attract the attention of many researchers. In the future, the examination of the properties of the numbers $y\left( n,\lambda \right) $ will continue and it will be investigated which other numbers and polynomials these numbers are related to.


\begin{thebibliography}{999}

	
\bibitem{Apostol} T. M. Apostol, \textit{On the Lerch zeta function}, Pacific J.
Math. \textbf{1}(2) (1951), 161--167.

\bibitem{balMS} Z. R. K. Balogh and M. J. Schlosser, \textit{$q$-Stirling numbers of
	the second kind $q$-Bell numbers for graphs}, Electronic Notes Discr. Math.
\textbf{54} (2016), 361--366.

\bibitem{Bayad} A. Bayad, \textit{Fourier expansions for Apostol-Bernoulli, Apostol-Euler and Apostol-Genocchi polynomials}, Math. Comp. \textbf{80}(276) (2011), 2219--2221.

\bibitem{Boyadv1} K. N. Boyadzhiev, \textit{A series transformation formula and
	related polynomials}, Int. J. Math. Math. Sci. \textbf{2005:23} (2005), 3849--3866.

\bibitem{Boyadv2} K. N. Boyadzhiev, \textit{Apostol-Bernoulli functions, derivative
	polynomials and Eulerian polynomials}, Adv. Appl. Discrete Math. \textbf{1}(2) (2008), 109--122.

\bibitem{Cangul} I. N. Cangul, H. Ozden, Y. Simsek, \textit{A new approach to $q$-Genocchi numbers and their interpolation functions}, Nonlinear Anal. \textbf{71}	(2009), e793--e799.

\bibitem{Carlitz} L. Carlitz, \textit{The number of derangements of a sequence with given specification}, Fibonacci. Quart. \textbf{16} (1978), 255--258.

\bibitem{Charalambos} C. A. Charalambides, \textit{Enumerative Combinatorics}, Chapman \& Hall/CRC, Boca Raton, London, New York, 2002.

\bibitem{CAC} C. A. Charalambides, \textit{Combinatorial Methods in Discrete
	Distributions}, Wiley-Interscience, Hoboken, New Jersey, 2005.

\bibitem{Choi1a} J. Choi, \textit{Explicit formulas for Bernoulli polynomials of order $n$}, Indian J. Pure Appl. Math. \textbf{27} (1996), 667--674.

\bibitem{ChoiJIA} J. Choi, \textit{Finite summation formulas involving binomial
	coefficients, harmonic numbers and generalized harmonic numbers}, J. Inequal. Appl. \textbf{2013:49} (2013), 1--11. 

\bibitem{ChoiSrivastavaTJM} J. Choi and H. M. Srivastava, \textit{The multiple Hurwitz zeta function and the multiple Hurwitz-Euler eta function}, Taiwanese J. Math. \textbf{15}(2) (2011), 501--522.

\bibitem{comtet} L. Comtet, \textit{Advanced Combinatorics}, D. Reidel Publication Company, Dordrecht-Holland/Boston-U.S.A., 1974.

\bibitem{Furdui} O. Furdui, \textit{Limits, Series, and Fractional Part Integrals, Problems in Mathematical Analysis}, Springer, 2013.

\bibitem{Guo} D. Guo and W. Chu, \textit{Summation formulae involving multiple Harmonic numbers}, Appl. Anal. Discrete Math. \textbf{15}(1) (2021), 201--212.

\bibitem{Dkim} D. S. Kim, T. Kim, J. Seo, \textit{A note on Changhee numbers and
	polynomials}, Adv. Stud. Theor. Phys. \textbf{7} (2013), 993--1003.

\bibitem{MSKim1a} M.-S. Kim, \textit{A note on sums of products of Bernoulli numbers}, Appl. Math. Lett. \textbf{24} (2011), 55--61.

\bibitem{KimSimsek2021} D. Kim, Y. Simsek, \textit{A New Family of Zeta Type
	Function Involving the Hurwitz Zeta Function and the Alternating Hurwitz
	Zeta Function}, Mathematics \textbf{9}(3) (2021), 233.

\bibitem{KimDahee} D. S. Kim, T. Kim, \textit{Daehee numbers and polynomials}, Appl. Math. Sci. \textbf{7}(120) (2013), 5969--5976.

\bibitem{T. Kim} T. Kim, \textit{$q$-Volkenborn Integration}, Russ. J. Math. Phys. \textbf{19} (2002), 288-299.

\bibitem{KIMjmaa2017} T. Kim, \textit{On the Analogs of Euler Numbers and Polynomials Associated with $p$-adic $q$-integral on $\mathbb{Z}_{\mathit{p}}$ at $q=1$}, J. Math. Anal. Appl. \textbf{331}(2) (2007), 779--792.

\bibitem{DSkimDaehee} D. S. Kim, T. Kim, \textit{Daehee Numbers	and Polynomials}, Appl. Math. Sci. (Ruse) \textbf{7}(120) (2013), 5969--5976.

\bibitem{Koepf} W. Koepf, \textit{Hypergeometric Summation, An Algorithmic Approach
	to Summation and Special Function Identities} (Second edition),
Springer-Verlag, London, 2014.

\bibitem{KucukogluAADM2019} I. Kucukoglu, Y. Simsek, \textit{On a family of special numbers and polynomials associated with Apostol-type numbers and
	poynomials and combinatorial numbers}, Appl. Anal. Discrete Math. \textbf{13} (2019), 478--494.

\bibitem{Ma} M. Ma, D. Lim, \textit{Degenerate Derangement Polynomials and Numbers}, Fractal Fract. \textbf{5}(3) (2021), 59.

\bibitem{GSimsek} G. V. Milovanovi\'{c}, Y. Simsek, V. S. Stojanovi\'{c}, \textit{A class of polynomials and connections with Bernoulli's numbers}, J. Anal. \textbf{27} (2019), 709--726.

\bibitem{ozdemirFilomat} G. Ozdemir, Y. Simsek, \textit{Generating functions for two-variable polynomials related to a family of Fibonacci type polynomials
	and numbers}, Filomat \textbf{30}(4) (2016), 969--975.

\bibitem{Ozdemir} G. Ozdemir, Y. Simsek, G. V. Milovanovi\'{c}, \textit{Generating
	functions for special polynomials and numbers including Apostol-type and Humbert-type polynomials}, Mediterr. J. Math. \textbf{14} (2017), 1--17.

\bibitem{Ozden} H. Ozden, I. N. Cangul, Y. Simsek, \textit{Multivariate interpolation functions of higher-order $q$-Euler numbers and their applications}, Abstr. Appl. Anal., Article ID 390857 (2008), 1--16.

\bibitem{Rao} K. S. Rao, V. Lakshminarayanan, \textit{Generalized Hypergeometric
	Functions Transformations and group theoretical aspects}, IOP Publishing Ltd., 2018.

\bibitem{RimJIAhar} S.-H. Rim, T. Kim, S.-S. Pyo, \textit{Identities between harmonic, hyperharmonic and Daehee numbers}, J. Inequal. Appl. \textbf{2018:168} (2018), 1--12.

\bibitem{5Riardon} J. Riordan, \textit{Introduction to Combinatorial Analysis}, Dover Publications, 2002.

\bibitem{Roman} S. Roman, \textit{The Umbral Calculus}, Dover Publications, New York, 2005.

\bibitem{Rota} G.-C. Rota, \textit{The number of partitions of a set}, American Math. Monthly \textbf{71}(5) (1964), 498--504. 

\bibitem{Schikof} W. H. Schikhof, \textit{Ultrametric Calculus: An Introduction to $p$-adic Analysis}, Cambridge Studies in Advanced Mathematics 4, Cambridge University Press, Cambridge, 1984.

\bibitem{SimsekJNT} Y. Simsek, \textit{$q$-Analogue of the twisted $l$-series and $q$-twisted Euler numbers}, J. Number Theory \textbf{110} (2005), 267--278.

\bibitem{Simsek JMAA} Y. Simsek, \textit{Twisted $(h,q)$-Bernoulli numbers and polynomials related to twisted $(h,q)$-zeta function and $L$-function}, J.
Math. Anal. Appl. \textbf{324}(2), 2006, 790--804.

\bibitem{Simsek1a} Y. Simsek, \textit{Multiple interpolation functions of higher order $(h;q)$-Bernoulli numbers}, AIP Conf. Proc. \textbf{1048} (2008), 486--489.

\bibitem{Simsek11a} Y. Simsek, \textit{$q$-Hardy Berndt type sums associated with $q$-Genocchi type zeta and $q$-$l$-functions}, Nonlinear Anal. \textbf{71} (2009), e377--e395.

\bibitem{SimsekJdea} Y. Simsek, \textit{Complete sum of products of $(h,q)$-extension of Euler polynomials and numbers}, J. Difference Equ. Appl. \textbf{16}(11) (2010), 1331--1348.

\bibitem{SimsekBook} Y. Simsek, \textit{Families of twisted Bernoulli numbers, twisted Bernoulli polynomials and their applications, Analytic Number
	Theory}, Approximation Theory, and Special Functions, Milovanovic, Gradimir
V.; Rassias, Michael Th. (Eds.), Eds., Springer, Berlin, pp. 149--214, 2014.

\bibitem{SimsekTJM2018} Y. Simsek, \textit{Construction of some new families of Apostol-type numbers and polynomials via Dirichlet character and $p$-adic $q$-integrals}, Turk. J. Math. \textbf{42} (2018), 557--577. 

\bibitem{simsekAADM} Y. Simsek, \textit{New families of special numbers for
	computing negative order Euler numbers and related numbers and polynomials}, Appl. Anal. Discrete Math. \textbf{12} (2018), 1--35.

\bibitem{Simsekfilomat} Y. Simsek, \textit{Combinatorial identities and sums for special numbers and polynomials}, Filomat \textbf{32}(20) (2018), 6869--6877.

\bibitem{simsekJMAA} Y. Simsek, \textit{Generating functions for finite sums involving higher powers of binomial coefficients: Analysis of hypergeometric
	functions including new families of polynomials and numbers}, J. Math. Anal. Appl. \textbf{477} (2019), 1328--1352.

\bibitem{simsekRSCM} Y. Simsek, \textit{Formulas for Poisson Charlier, Hermite, Milne-Thomson and other type polynomials by their generating functions and $p$-adic integral approach}, Rev. R. Acad. Cienc. Exactas F\'{i}s. Nat. Ser. A Mat. RACSAM \textbf{113}(2) (2019), 931--948.

\bibitem{SimsekMTJPM} Y. Simsek, \textit{Explicit formulas for $p$-adic integrals: Approach to $p$-adic distributions and some families of special
	numbers and polynomials}, Montes Taurus J. Pure Appl. Math. \textbf{1}(1) (2019), 1--76.

\bibitem{SimsekMTJPAM2020} Y. Simsek, \textit{Interpolation functions for new
	classes special numbers and polynomials via applications of $p$-adic	integrals and derivative operator}, Montes Taurus J. Pure Appl. Math. \textbf{3}(1) (2021), 1--24.

\bibitem{SimsekREVISTA} Y. Simsek, \textit{New integral formulas and identities involving special numbers and functions derived from certain class of
	special combinatorial sums}, Rev. R. Acad. Cienc. Exactas F\'{i}s. Nat. Ser. A Mat. RACSAM \textbf{115:66} (2021), 1--14.

\bibitem{SimsekBull2021} Y. Simsek, \textit{Miscellaneous formulae for the certain class of combinatorial sums and special numbers}, Bull. Cl. Sci. Math. Nat. Sci. Math. \textbf{45} (2020).

\bibitem{SimsekASCM2021} Y. Simsek, \textit{Construction of generalization Leibnitz type numbers and their properties}, Adv. Stud. Contemp. Math. (Kyungshang) \textbf{31}(3) (2021), 311--323.

\bibitem{1Sofo} A. Sofo, \textit{Quadratic alternating harmonic number sums}, J. Number Theory \textbf{154} (2015), 144--159.

\bibitem{SrivatavaChoi} H. M. Srivastava, J. Choi, \textit{Zeta and $q$-Zeta
	Functions and Associated Series and Integrals}, Elsevier Science Publishers,
Amsterdam, 2012.	

\bibitem{RT} R. Tremblay, B. J. Fugere, \textit{Products of two restricted
	hypergeometric functions}, J. Math. Anal. Appl. \textbf{198}(3) (1996), 844--852.

\bibitem{ChoiMTJPAM} T. Usman, N. Khan, M. Saif, J. Choi, \textit{A Unified Family of Apostol-Bernoulli Based Poly-Daehee Polynomials}, Montes Taurus J.
Pure Appl. Math. \textbf{3}(3) (2021), 1--11.

\bibitem{OEIS} \texttt{https://oeis.org/A025529}

\bibitem{Volkenborn} A. Volkenborn, \textit{On Generalized $p$-adic Integration}, M\'{e}m. Soc. Math. Fr. \textbf{39-40} (1974), 375--384.

\bibitem{Zave} D. A. Zave, \textit{A series expansion involving the harmonic numbers}, Inform. Process. Lett. \textbf{5}(3) (1976), 75--77.

\end{thebibliography}
\end{document}